\documentclass[aop,MSNbibl,seceqn,dvips]{arximspdf}
\usepackage{mathrsfs}

\doi{10.1214/14-AOP922} 
\volume{42}
\issue{5}
\pubyear{2014}
\firstpage{1952}
\lastpage{1979}

\makeatletter

\def\idnk{\operatorname{id}_{k,n}}
\def\idk{\operatorname{id}_k}
\def\canonk{\mathbf{Z}_k}
\def\canonkn{\mathbf{Z}_{k,n}}
\def\Rn{\mathbf{R}_n}
\def\Rmn{\mathbf{R}_{m,n}}
\def\Frag{\textsc{Frag}}
\def\Coag{\textsc{Coag}}
\def\cutpaste{\textsc{Cut-Paste}}
\def\PD{\operatorname{PD}}
\def\partitionsNk{\mathcal{P}_{\mathbb{N}:k}}
\def\partitionsN{\mathcal{P}_{\mathbb{N}}}
\def\symmetrick{\mathscr{S}_k}

\def\symmetricN{\mathscr{S}_{\mathbb{N}}}

\def\stochk{\mathcal{S}_k}
\def\simplexk{\Delta_k}
\def\equalinlaw{=_{\mathcal{L}}}
\def\Nb{\mathbb{N}}
\def\Bb{\mathcal{B}}
\def\Pb{\mathbf{P}}
\def\Pbb{\mathbb{P}}
\def\Xb{\mathbf{X}}
\def\Pib{\bolds{\Pi}}
\def\Mb{\mathbf{M}}

\newtheorem{thmm}{Theorem}[section]
\newtheorem{lemma}[thmm]{Lemma}
\newtheorem{prop}[thmm]{Proposition}
\newtheorem{cor}[thmm]{Corollary}
\newproclaim{defn}[thmm]{Definition}
\newproclaim{example}[thmm]{Example}

\newproclaim{rmk}[thmm]{Remark}
\newcommand{\eqref}[1]{(\ref{#1})}
\makeatother

\begin{document}
\begin{frontmatter}

\title{The cut-and-paste process\thanksref{T1}}
\runtitle{The cut-and-paste process}
\thankstext{T1}{Supported in part by NSF Grant DMS-13-08899 and NSA
Grant H98230-13-1-0299.}

\begin{aug}
\author{\fnms{Harry} \snm{Crane}\corref{}\ead[label=e1]{hcrane@stat.rutgers.edu}}
\runauthor{H. Crane}
\affiliation{Rutgers University}
\address{Department of Statistics\\
Rutgers University\\
110 Frelinghuysen Road\\
Piscataway, New Jersey 08854\\
USA\\
\printead{e1}}
\end{aug}

\received{\smonth{7} \syear{2012}}
\revised{\smonth{1} \syear{2014}}

%
\begin{abstract}
We characterize the class of exchangeable Feller processes evolving on
partitions with boundedly many blocks.
In continuous-time, the jump measure decomposes into two parts: a
$\sigma$-finite measure on stochastic matrices and a collection of
nonnegative real constants.
This decomposition prompts a L\'evy--It\^o representation.
In discrete-time, the evolution is described more simply by a product
of independent, identically distributed random matrices.
\end{abstract}

%

\pdfkeywords{60J25, 60G09, 60J35, Exchangeable random partition, de Finetti's theorem, Levy-Ito decomposition, paintbox process, coalescent process, interacting particle system, Feller process, random matrix product}
\end{frontmatter}

\section{Introduction}\label{sec:intro}

For fixed $k=1,2,\ldots,$ a \emph{$k$-coloring} of $\Nb:=\{
1,2,\ldots
\}$
is an infinite sequence $x=x^1x^2\cdots$ taking values in $[k]:=\{
1,\ldots,k\}$.
Two operations bear on our main theorems:
\begin{itemize}
\item relabeling: for any permutation $\sigma\dvtx\Nb\rightarrow\Nb
$, the
\emph{relabeling} of $x=x^1x^2\cdots$ by $\sigma$ is
%
%
\begin{equation}
\label{eq:relabeling} x^{\sigma}:=x^{\sigma(1)}x^{\sigma(2)}\cdots\quad
\mbox{and}
\end{equation}
\item restriction: for any finite $n=1,2,\ldots,$ the \emph{restriction}
of $x$ to a $k$-coloring of $[n]$ is
%
%
\begin{equation}
\label{eq:restriction} x^{[n]}:=x^1\cdots x^n.
\end{equation}
\end{itemize}
A Markov process $\Xb=(X_t, t\geq0)$ on $[k]^{\mathbb{N}}$, the
space of
infinite $k$-colorings, is
\begin{longlist}[(A)]
\item[(A)] \emph{exchangeable} if $\Xb^{\sigma}=(X_t^{\sigma},
t\geq0)$
is a version of $\Xb$ for all finite permutations $\sigma\dvtx\Nb
\rightarrow
\Nb$ and
\item[(B)] \emph{consistent (under subsampling)} if $\Xb
^{[n]}=(X_t^{[n]}, t\geq0)$ is a Markov chain on $k$-colorings of
$[n]$, for all finite $n=1,2,\ldots.$
\end{longlist}

We characterize both $[k]^{\mathbb{N}}$-valued Markov processes satisfying
(A) and (B) and a class of partition-valued processes with analogous properties.
When $[k]^{\mathbb{N}}$ is endowed with the product-discrete topology,
exchangeability and consistency are equivalent to
exchangeability and the Feller property;
and so our main theorems characterize exchangeable Feller processes on
$[k]^{\mathbb{N}}$ and $\partitionsNk$, partitions of $\Nb$ with at most
$k$ blocks.

\subsection{Discrete-time characterization}
A stochastic matrix $S=(S_{\mathit{ii}'}, 1\leq i,i'\leq k)$ has nonnegative
entries and all rows summing to one, and it determines the transition
probabilities of a time-homogeneous Markov chain $\mathbf{Y}=(Y_m,
m\geq0)$ on $[k]$ by
%
%
\begin{equation}
\label{eq:tps-S} \Pbb_S\bigl\{Y_1=i'\mid
Y_0=i\bigr\}=S_{\mathit{ii}'},\qquad i,i'=1,\ldots,k.
\end{equation}
From any probability measure $\Sigma$ on the space of $k\times k$
stochastic matrices, we construct a Markov chain $\Xb^*_{\Sigma
}:=(X_m^*, m\geq0)$ on $[k]^{\mathbb{N}}$ as follows.
First, we let $X_0^*$ be an exchangeable initial state and
$S_1,S_2,\ldots$ be independent, identically distributed (i.i.d.)
random matrices from $\Sigma$.
Then, for $m=1,2,\ldots,$ we generate the components of
$X_m^*=X_{m}^{*1}X_{m}^{*2}\cdots,$ given $X_{m-1}^*,\ldots,X_0^*,S_1,S_2,\ldots,$ conditionally independently from transition
probability matrix $S_{m}:=(S_m(i,i'), 1\leq i,i'\leq k)$,
\[
\mathbb{P}\bigl\{X_{m}^{*j}=i'\mid
X_{m-1}^*, S_{m}\bigr\}=S_m\bigl(i,i'
\bigr) \qquad\mbox {on the event }X_{m-1}^{*j}=i.
\]
Such a construction exists for all exchangeable and consistent Markov
chains on~$[k]^{\mathbb{N}}$.

%
\begin{thmm}\label{thmm:cut-and-paste chain}
Let $\Xb=(X_m, m\geq0)$ be a discrete-time, exchangeable, consistent
Markov chain on $[k]^{\mathbb{N}}$. Then there exists a unique probability
measure $\Sigma$ such that $\mathbf{X}^*_{\Sigma}$ is a version of
$\Xb$.
\end{thmm}

To any $x\in[k]^{\mathbb{N}}$, the \emph{asymptotic frequency vector}
$|x|:=(f_1(x),\ldots,f_k(x))$ is an element of the $(k-1)$-dimensional
simplex $\simplexk$, where
%
%
\begin{equation}
\label{eq:asymptotic frequency} f_i(x):=\lim_{n\rightarrow\infty}n^{-1}
\sum_{j=1}^n\mathbf{1}\bigl\{
x^j=i\bigr\},\qquad  i=1,\ldots,k,
\end{equation}
is the limiting proportion of coordinates labeled $i$ in $x$, if it exists.
With probability one, the asymptotic frequency vector of any
exchangeable $k$-coloring exists and $|\Xb|:=(|X_m|, m\geq0)$ is a
sequence in $\simplexk$.
From the same i.i.d. sequence $S_1,S_2,\ldots$ used to construct $\Xb
_{\Sigma}^*=(X_m^*, m\geq0)$ in Theorem~\ref{thmm:cut-and-paste
chain}, we can construct $\bolds{\Phi}_{\Sigma}:=(\Phi_m, m\geq
0)$ in
$\simplexk$ by putting $\Phi_0:=|X_0^*|$ and
%
%
\begin{equation}
\label{eq:Phi} \Phi_m:=\Phi_{m-1}S_m=
\Phi_0S_1\cdots S_m,\qquad m\geq1,
\end{equation}
where $\Phi_{m-1}S_m$ in \eqref{eq:Phi} is the usual right action of a
$k\times k$ matrix on a $1\times k$ row vector.

%
\begin{thmm}\label{thmm:induced simplex chain}
Let $\Xb=(X_m, m\geq0)$ be a discrete-time, exchangeable, consistent
Markov chain on $[k]^{\mathbb{N}}$.
Then $\bolds{\Phi}_{\Sigma}$ is a version of $|\Xb|$, where
$\Sigma$
is the unique probability measure from Theorem~\ref{thmm:cut-and-paste chain}.
\end{thmm}

Together, Theorems \ref{thmm:cut-and-paste chain} and \ref{thmm:induced
simplex chain} relate the evolution of discrete-time Markov chains to
products of i.i.d. random matrices.
Crane and Lalley \cite{CraneLalley2012a} have combined representation
\eqref{eq:Phi} with the Furstenberg--Kesten theorem \cite
{FurstenbergKesten1960} to identify a class of these chains that
exhibits the cutoff phenomenon.

\subsection{Continuous-time characterization}\label{section:intro:Levy-Ito}

In continuous-time, an exchangeable, consistent Markov process $\Xb
=(X_t, t\geq0)$ can jump infinitely often, and thus, behaves
differently than its discrete-time counterpart; but consistency limits
this behavior: since each restriction $\Xb^{[n]}$ is a finite state
space Markov process, it must remain in each visited state for a
positive amount of time.
The upshot of these observations is a characterization of the
transition law of $\Xb$ by a unique $\sigma$-finite measure on
$k\times
k$ stochastic matrices and a unique collection of nonnegative constants.

Our next theorem yields a L\'evy--It\^o-type characterization of $\Xb$
by dividing its discontinuities into two cases.
Let $t>0$ be the time of a discontinuity in $\Xb$.
Then either
\begin{longlist}[(II)]
\item[(I)] a positive proportion of coordinates changes colors at time
$t$, that is, %
\[
\lim_{n\rightarrow\infty}n^{-1}\sum_{j=1}^n
\mathbf{1}\bigl\{ X_{t-}^j\neq X_{t}^{j}
\bigr\}>0 \quad\mbox{or}
\]
\item[(II)] a zero proportion of coordinates changes colors at time
$t$, that is, %
\[
\lim_{n\rightarrow\infty}n^{-1}\sum_{j=1}^n
\mathbf{1}\bigl\{ X_{t-}^j\neq X_{t}^{j}
\bigr\}=0.
\]
\end{longlist}
In discrete-time, Type-(I) jumps are governed by a probability measure
$\Sigma$ and Type-(II) transitions are forbidden.
In continuous-time, Type-(I) jumps are governed by a $\sigma$-finite
measure $\Sigma$ and Type-(II) transitions include only \emph
{single-index flips}, that is, jumps for which exactly one coordinate
changes color.
Deciding the Type-(II) jump rates is a collection of nonnegative
constants $\mathbf{c}=(\mathbf{c}_{\mathit{ii}'}, 1\leq i\neq i'\leq k)$:
independently, each coordinate changes colors from $i$ to $i'$ at rate
$\mathbf{c}_{\mathit{ii}'}$.
The transition law of $\Xb$ is characterized by the pair $(\Sigma,\mathbf{c})$.

We do not fully explain $(\Sigma,\mathbf{c})$ and its relation to
$\Xb$
until Section~\ref{section:Levy-Ito}.
Sparing the details, we write $\Xb_{\Sigma,\mathbf{c}}^*$ to denote a
continuous-time Markov process constructed from a Poisson point process
with intensity measure determined by $(\Sigma,\mathbf{c})$.
Theorem~\ref{thmm:Levy-Ito} says that any exchangeable, consistent
Markov process $\Xb$ admits a version with this construction.

%
\begin{thmm}\label{thmm:Levy-Ito}
Let $\Xb=(X_t, t\geq0)$ be a continuous-time, exchangeable, consistent
Markov process on $[k]^{\mathbb{N}}$.
Then there exists a unique measure $\Sigma$ satisfying~\eqref
{eq:regularity Sigma}
and unique nonnegative constants $\mathbf{c}=(\mathbf{c}_{\mathit{ii}'}, 1\leq
i\neq i'\leq k)$ such that $\Xb_{\Sigma,\mathbf{c}}^*$ is a version of
$\Xb$.
\end{thmm}

In Theorem~\ref{thmm:Levy-Ito}, $\Sigma$ is required to satisfy
%
%
\begin{equation}
\label{eq:regularity Sigma} \Sigma\bigl(\{I_k\}\bigr)=0 \quad\mbox{and}\quad \int
_{\stochk}(1-S_*)\Sigma (dS)<\infty,
\end{equation}
where $I_k$ is the $k\times k$ identity matrix, $S_*:=\min
(S_{11},\ldots,S_{kk})$ for any $k\times k$ stochastic matrix $S$, and $\stochk$ is
the space of $k\times k$ stochastic matrices.
Consistency imposes~\eqref{eq:regularity Sigma}: uniqueness requires
the first half, finiteness of finite-dimensional jump rates forces the
second half.

As in discrete-time, we define the projection of $\Xb=(X_t, t\geq0)$
into $\simplexk$ by $|\Xb|=(|X_t|, t\geq0)$.
Unlike discrete-time, the existence of $|\Xb|$ does not follow directly
from de Finetti's theorem because now $\Xb$ is an uncountable collection.

%
\begin{thmm}\label{thmm:induced simplex process}
Let $\Xb=(X_t, t\geq0)$ be a continuous-time, exchangeable, consistent
Markov process on $[k]^{\mathbb{N}}$.
Then $|\Xb|=(|X_t|, t\geq0)$ exists almost surely and is a Feller
process on $\simplexk$.
\end{thmm}

Theorems \ref{thmm:Levy-Ito} and \ref{thmm:induced simplex process}
give the L\'evy--It\^o representation.
The projection $|\Xb|$ jumps only at the times of Type-(I)
discontinuities in $\Xb$; at other times, it follows a continuous,
deterministic trajectory.
Thus, Theorem~\ref{thmm:Levy-Ito} warrants the heuristic interpretation
that $\Sigma$ governs the ``discrete'' component of $\Xb$ and
$\mathbf
{c}$ governs the ``continuous'' component.

\subsection{Partition-valued Markov processes}
Any $x\in[k]^{\mathbb{N}}$ determines a partition $\pi=\Bb(x)$ of
$\Nb$ through
%
%
\begin{equation}
\label{eq:coloring to partition} i\mbox{ and }j\mbox{ are in the same block of }\pi
\quad\Longleftrightarrow\quad x^i=x^j.
\end{equation}
If the characteristic pair $(\Sigma,\mathbf{c})$ treats colors
symmetrically, that is, $\Sigma$ is row--column exchangeable and
$\mathbf{c}_{\mathit{ii}'}=\mathbf{c}_{jj'}=c$ for all $i\neq i'$ and $j\neq
j'$, then the projection $\Bb(\Xb_{\Sigma,\mathbf{c}}^*)=(\Bb(X_t^*),
t\geq0)$ into $\partitionsNk$ through \eqref{eq:coloring to partition}
is an exchangeable, consistent Markov process on $\partitionsNk$. Our
main theorem for partition-valued processes states that any
exchangeable, consistent Markov process on $\partitionsNk$ can be
generated by projecting an exchangeable, consistent Markov process from
$[k]^{\mathbb{N}}$.

%
\begin{thmm}\label{thmm:Levy-Ito partition}
Let $\bolds{\Pi}$ be an exchangeable, consistent Markov process on~$\partitionsNk$.
\begin{itemize}
\item In discrete-time, there exists a unique, row--column exchangeable
probability measure $\Sigma$ such that $\Bb(\Xb^*_{\Sigma})$ is a
version of $\bolds{\Pi}$;
\item in continuous-time, there exists a unique, row--column
exchangeable measure satisfying \eqref{eq:regularity Sigma} and a
unique constant $c\geq0$ such that $\Bb(\Xb_{\Sigma,\mathbf
{c}}^*)$ is
a version of~$\bolds{\Pi}$, where $\mathbf{c}_{\mathit{ii}'}=c$ for all
$1\leq
i\neq i'\leq k$.
\end{itemize}
\end{thmm}

Analogously to \eqref{eq:asymptotic frequency}, we define the
asymptotic frequency of $\pi\in\partitionsN$ by $|\pi|^{\downarrow}$,
the asymptotic block frequencies of $\pi$ in decreasing order of size.
When it exists, $|\pi|^{\downarrow}$ is an element of the \emph{ranked
$k$-simplex} $\Delta_k^{\downarrow}$.

%
\begin{thmm}\label{thmm:induced simplex process-partition}
Let $\bolds{\Pi}=(\Pi_t, t\geq0)$ be a continuous-time, exchangeable,
consistent Markov process on $\partitionsNk$.
Then $|\bolds{\Pi}|^{\downarrow}:=(|\Pi_t|^{\downarrow}, t\geq0)$
exists almost surely and is a Feller process on $\Delta_k^{\downarrow}$.
\end{thmm}

\subsection{The cut-and-paste process}
We call $\Xb^*_{\Sigma,\mathbf{c}}$ a \emph{cut-and-paste process}: its
jumps occur by first \emph{cutting} each color class into subclasses and
then \emph{pasting} subclasses together.
When $(\Sigma,\mathbf{c})$ treats colors symmetrically, we call $\Xb
^*_{\Sigma,\mathbf{c}}$ and its projection into $\partitionsNk$ a
\emph{homogeneous cut-and-paste process}.

Cut-and-paste processes should not be conflated with synonymous, but
not analogous, \emph{split-and-merge} \cite{Pitman2002a} and \emph
{coagulation--fragmentation} processes \cite{DurrettGranovskyGueron1999}.
The latter processes share aspects, but are not one, with the
cut-and-paste process.
Each process evolves by operations that divide (cut, split, fragment)
and unite (paste, merge, coagulate), but split-and-merge processes
evolve on interval partitions, coagulation--fragmentation processes on
set partitions, and cut-and-paste processes on $k$-colorings.
At the time of a jump, a cut-and-paste process undergoes two operations
simultaneously (cut \emph{and} paste), the others undergo only one
operation (split \emph{or} merge, coagulate \emph{or} fragment).

Theorems \ref{thmm:Levy-Ito partition} and \ref{thmm:induced simplex
process-partition} do elicit qualitative connections to exchangeable
coalescent and fragmentation processes \cite{Bertoin2001a,Pitman1999a},
both of which are characterized by pairs $(\nu,c)$, where $\nu$ is a
unique $\sigma$-finite measure on ranked-mass partitions and $c\geq0$
is a unique constant.
For coalescent processes, $\nu$ determines the rate of multiple
collisions and $c$ the rate of binary coalescence.
For fragmentation processes, $\nu$ determines the rate of dislocation
and $c$ the rate of erosion.
In both cases, $(\nu,c)$ gives a L\'evy--It\^o description.
But, in a strict sense, processes on $\partitionsNk$ behave differently
than those on $\partitionsN$ \cite{Crane2011a,Crane2012c}, and
Theorem~\ref{thmm:Levy-Ito partition} neither refines nor is a special
case of
previous results.
In Section~\ref{section:relation to cfp}, we further discuss any
relationships (and lack thereof) between cut-and-paste, coalescent and
fragmentation processes.

\subsection{Applications to DNA sequencing}\label{section:DNA sequencing}
Decades ago, population genetics applications motivated the initial
study of random partitions and partition-valued processes \cite
{Ewens1972,Kingman1978a,Kingman1982}.
Somewhat later, Bertoin \cite{Bertoin2001a,Bertoin2006} and Pitman
\cite{Pitman1999a,Pitman2005} connected coalescent and fragmentation
processes to Brownian motion, L\'evy processes and subordinators.
In the present, DNA sequencing inspires processes restricted to
partitions with a bounded number of blocks.

For let the colors correspond to DNA nucleotides, adenine (A), cytosine
(C), guanine (G) and thymine (T).
Then, for a sample of $n$ individuals, $X^1\cdots X^n\in\{A,C,G,T\}
^{[n]}$ is a string of DNA nucleotides at a particular chromosomal
site, where $X^i$ denotes the nucleotide of individual $i=1,\ldots,n$.
If we observe a DNA sequence $(X^i_m, m\geq0)$ for each $i=1,\ldots,n$, then $(X_m, m\geq0)$ is a sequence in $\{A,C,G,T\}^{[n]}$, with
$X_m=X_m^1\cdots X_m^n$.
By forgetting colors (in this case nucleotides), we obtain a sequence
of set partitions; see Table~\ref{table:0}.

In practice, biological phenomena such as recombination induce
dependence among nearby chromosomal sites.
For modeling this dependence, the Markov property strikes a balance
between practical feasibility and mathematical tractability.
Exchangeability and consistency incorporate a logical structure that is
apt for DNA sequencing.
See \cite{Crane2014a} for a detailed statistical consideration of these
applications.

%
\begin{table}
\caption{An array of DNA sequences for 3 individuals. From this array,
we obtain a sequence in $\{A,C,G,T\}^{[3]}$: (AAT, ATT, TTT, CCG, CGG,
GGC, AAT, $\ldots$).
By ignoring nucleotide labels, we obtain the sequence $(12|3, 1|23,
123, 12|3, 1|23, 12|3, 12|3,\ldots)$ of partitions of the set $\{
1,2,3\}$}\label{table:0}
\begin{tabular*}{\textwidth}{@{\extracolsep{\fill}}lcccccccc@{}}
\hline
\textbf{Individuals/sites}& $\bolds{m=1}$ & \textbf{2} & \textbf{3} & \textbf{4} &
\textbf{5} & \textbf{6} & \textbf{7} & $\bolds{\cdots}$\\
\hline
$X_m^1$ & A & A & T & C & C & G & A & $\cdots$ \\[2pt]
$X_m^2$ & A & T & T & C & G & G & A & $\cdots$ \\[2pt]
$X_m^3$ & T & T & T & G & G & C & T & $\cdots$\\
\hline
\end{tabular*}
\end{table}

\subsection{Discussion of main theorems}\label{section:discussion}
For concreteness, let $\Xb$ be a discrete-time Markov chain on $\{1,2\}
^{\Nb}$.
According to Theorem~\ref{thmm:cut-and-paste chain}, a transition
$X\mapsto X'$ can be generated in two steps:
\begin{enumerate}[(ii)]
\item[(i)] Draw a random pair $(p_1,p_2)$ of success probabilities from
a probability measure $\Sigma$ on $[0,1]\times[0,1]$.
\item[(ii)] Given $(p_1,p_2)$, update each coordinate $j=1,2,\ldots$ of
$X$ independently by the following coin flipping process.
\begin{itemize}[$-$]
\item[$-$] If $X^j=1$, flip a $p_1$-coin ($\Pbb\{\mbox{heads}\}=p_1)$;
otherwise, flip a $p_2$-coin.
\item[$-$] If the outcome is heads, put $X'^j=1$; otherwise, put $X'^j=2$.
\end{itemize}
\end{enumerate}
The pair $(p_1,p_2)$ determines a $2\times2$ stochastic matrix
\[
S= %
\pmatrix{ p_1 & 1-p_1\vspace*{2pt}
\cr
p_2 & 1-p_2 },
\]
which describes the transition probability matrix for each coordinate,
as in \eqref{eq:tps-S}.
By the law of large numbers, the proportion of coordinates labeled 1 in
$X'$ equals
\begin{eqnarray*}
f_1\bigl(X'\bigr)&=&\Pbb\bigl\{X'^1=1
\mid X^1=1\bigr\}f_1(X)+\Pbb\bigl\{X'^1=1
\mid X^1=2\bigr\} f_2(X)\\
&=&p_1f_1(X)+p_2f_2(X).
\end{eqnarray*}
Overall, the asymptotic frequencies $|X'|=(f_1(X'),f_2(X'))$ of $X'$
are the entries of
\[
|X|S= %
\pmatrix{ f_1(X) & f_2(X)} %
\pmatrix{ p_1 & 1-p_1\vspace*{2pt}
\cr
p_2 &
1-p_2 }.
\]

In discrete-time, exchangeability implies that if $X'\neq X$, then the
proportion of coordinates changing colors from $X$ to $X'$ is strictly positive.
In continuous-time, the transition rate $X\mapsto X'$ need not be
bounded, and thus, $\Sigma$ need not be finite.
Furthermore, there is no requirement that a strictly positive
proportion of coordinates changes colors at the time of a discontinuity.
However, the consistency assumption implies that any finite collection
of coordinates jumps at a finite rate, producing condition \eqref
{eq:regularity Sigma}.
Together, exchangeability and consistency restrict Type-(II)
discontinuities to involve only a single coordinate, called a \emph
{single-index flip}.
For instance, if ``double-index flips'' were permitted, that is, a
pair of indices changes colors simultaneously while all other
coordinates remain unchanged, then the finite restrictions of $\mathbf
{X}$ could not be c\`adl\`ag.
To see this, suppose any pair $(X^n,X^{n'})$, $n<n'$, changes from
$(1,1)$ to $(2,2)$ at positive rate $\mathbf{r}$.
Then, by exchangeability, any pair $(X^n,X^{n'+j})$, $j\geq1$, in state
$(1,1)$ must also flip at rate $\mathbf{r}$.
For any such jump, the restriction of $\Xb$ to $[n]$ witnesses only a
change in coordinate $n$ at rate $\sum_{n'> n}\mathbf{r}=\infty$, which
contradicts assumption (B).
For similar reasons, condition \eqref{eq:regularity Sigma} prevents
infinitely many Type-(I) discontinuities from bunching up in any finite
restriction of $\Xb$.

Upon observing our main theorems for $[k]^{\mathbb{N}}$-valued
processes, the
analogous conclusions for $\partitionsNk$-valued processes are nearly
immediate.
The key observation is that the projection of $\mathbf{X}$ into
$\partitionsNk$ preserves the Markov property only if the transition
law of $\mathbf{X}$ treats the labels $[k]$ symmetrically, which
requires row--column exchangeability of $\Sigma$ and $\mathbf
{c}_{\mathit{ii}'}=\mathbf{c}_{jj'}=c$ for all $i\neq i'$, $j\neq j'$.

\subsection{Examples}\label{section:examples}
We illustrate our main theorems with three examples: two exchangeable,
consistent Markov processes on $[k]^{\mathbb{N}}$ (one in
discrete-time and
one in continuous-time) and a family of exchangeable Markov chains that
is not consistent (the Ehrenfest walk on the hypercube).
Example~\ref{ex:non-example} shows why discrete-time chains cannot
admit single-index flips.

%
\begin{example}[(A reversible discrete-time chain \cite
{Crane2011a})]\label{ex:discrete}
For $\alpha>0$, we define transition probabilities
%
%
\begin{equation}
\label{eq:example tps} P_n\bigl(x,x'\bigr):=\prod
_{i=1}^k\frac{\prod_{i'=1}^k(\alpha/k)^{\uparrow
\mathbf
{n}_{\mathit{ii}'}(x,x')}}{\alpha^{\uparrow\mathbf{n}_i(x)}},\qquad x,x'\in
\mathop{[k]^{[n]}},
\end{equation}
where $\mathbf{n}_{\mathit{ii}'}(x,x'):=\#\{j\in[n]\dvtx x^j=i\mbox{ and }
x'^j=i'\}
$, $\mathbf{n}_i(x):=\#\{j\in[n]\dvtx x^j=i\}$, and $\alpha
^{\uparrow
n}:=\alpha(\alpha+1)\cdots(\alpha+n-1)$.
This transition probability is reversible with respect to
\[
\lambda_{\xi}^{(n)}(x)=\frac{\prod_{i=1}^{k}\alpha^{\uparrow
\mathbf
{n}_i(x)}}{(k\alpha)^{\uparrow n}},\qquad x\in
\mathop{[k]^{[n]}},
\]
and projects to a transition probability on $\mathop{\mathcal
{P}_{[n]:k}} $ (partitions
of $[n]$ with at most $k$ blocks) with reversible stationary distribution
\[
\varrho_{\xi}^{(n)}(\pi):=\frac{k!}{(k-\#\pi)!}\frac{\prod_{b\in
\pi
}\alpha^{\uparrow\#b}}{(k\alpha)^{\uparrow n}},\qquad
\pi\in\mathop {\mathcal{P}_{[n]:k}},
\]
where $\#\pi$ denotes the number of blocks of $\pi$ and $\#b$ denotes
the cardinality of $b\subseteq[n]$.

Namely, in Theorem~\ref{thmm:cut-and-paste chain}, the transition
probabilities in \eqref{eq:example tps} correspond to the homogeneous
cut-and-paste chain with $\Sigma_{\alpha/k}=\xi_{\alpha/k}\otimes
\cdots
\otimes\xi_{\alpha/k}$, where $\xi_{\alpha}$ is the symmetric Dirichlet
distribution with parameter $(\alpha,\ldots,\alpha)$.
That is, $S\sim\Sigma_{\alpha/k}$ is a random matrix whose rows are
independent and identically distributed from $\operatorname
{Dirichlet}(\alpha/k,\ldots,\alpha/k)$.
\end{example}

%
\begin{example}[(A purely continuous process)]\label{ex:continuous}
For $\mathbf{c}_{12},\mathbf{c}_{21}>0$, let each coordinate of $\Xb
=(X_t, t\geq0)$ evolve independently, each jumping from 1 to 2 at rate
$\mathbf{c}_{12}$ and from 2 to 1 at rate $\mathbf{c}_{21}$.
The projection of $\Xb$ into the simplex evolves continuously and
deterministically by a constant interchange of mass between the colors
1 and 2.
Eventually, the projection settles to the fixed point
\[
\biggl(\frac{c_{21}}{c_{12}+c_{21}},\frac
{c_{12}}{c_{12}+c_{21}} \biggr).
\]

The projection into $\partitionsNk$ is Markov only if $c_{12}=c_{21}$.
In this case, the projection settles to $(1/2,1/2)$ and, in
equilibrium, there is a constant and equal flow of mass between the two blocks.
\end{example}

%
\begin{example}[(Nonexample: Ehrenfest chain on $\{0,1\}^{[n]}$)]\label
{ex:non-example}
The family of discrete-time Ehrenfest chains on the hypercubes $\{0,1\}
^{[n]}$, $n\in\mathbb{N}$, is not consistent, and thus, not covered by
our theory.
On $\{0,1\}^{[n]}$, an Ehrenfest chain $\mathbf{X}^{[n]}$ evolves by
choosing a coordinate $1,\ldots,n$ uniformly at random and then
flipping a fair coin to decide its value at the next time. All other
coordinates remain unchanged. In the language of Section~\ref
{section:discussion}, all transitions of this chain are single-index flips.

The finite-dimensional chains are exchangeable but not consistent.
For any \mbox{$n\in\Nb$}, the probability that $\mathbf{X}^{[n]}$ remains in
the same state after a transition is $1/2$, whereas the projection of
an Ehrenfest chain $\mathbf{X}^{[n+1]}$ on $\{0,1\}^{[n+1]}$ into $\{
0,1\}^{[n]}$ remains in the same state with probability
$(n+2)/(2n+2)\neq1/2$.
\end{example}

Six sections compose the paper. In Section~\ref{section:preliminaries},
we lay out definitions and notation; in Section~\ref{section:random
matrix product}, we establish Theorems \ref{thmm:cut-and-paste chain}
and \ref{thmm:induced simplex chain}; in Section~\ref
{section:Levy-Ito}, we prove Theorems \ref{thmm:Levy-Ito} and \ref
{thmm:induced simplex process}; in Section~\ref{section:partition
chains}, we deduce Theorems \ref{thmm:Levy-Ito partition} and \ref
{thmm:induced simplex process-partition}; in Section~\ref
{section:conclusion}, we conclude.

\section{Preliminaries}\label{section:preliminaries}
\subsection{Notation}\label{section:notation}
Throughout the paper, we write $x$ to denote a $k$-coloring, $X$ a
random $k$-coloring and $\mathbf{X}$ a random collection of $k$-colorings.
We write $\pi$ to denote a partition, $\Pi$ a random partition, and
$\bolds{\Pi}$ a random collection of partitions.
For terminology and notation pertaining to both $k$-colorings and
partitions, we write $\lambda$, $\Lambda$, and $\bolds{\Lambda}$, as
appropriate.
A collection $\bolds{\Lambda}=(\Lambda_m, m\geq0)$ indexed by $m$
evolves in discrete-time, that is, $m=1,2,\ldots,$ and $\bolds
{\Lambda
}=(\Lambda_t, t\geq0)$ indexed by $t$ evolves in continuous-time, that
is, $t\in[0,\infty)$.

\subsection{Partitions and colorings}\label{section:set partitions}
For fixed $k\in\mathbb{N}$, a \emph{$k$-coloring} of $[n]=\{1,\ldots,n\}
$ is a $[k]$-valued sequence $x=x^1\cdots x^n$.
A \emph{partition} of $[n]$ is a collection $\pi=\{B_1,\ldots,B_r\}$ of
nonempty, disjoint subsets (blocks) satisfying $\bigcup_{j=1}^r B_j=[n]$.
We can also regard $\pi$ as an equivalence relation $\sim_\pi$, where
\[
i\sim_\pi j \quad\Longleftrightarrow\quad i\mbox{ and }j\mbox{ are in the same
block of }\pi.
\]
Upon removal of its colors, any $k$-coloring $x$ projects to a unique
partition $\Bb_n(x)$ of $[n]$, as in \eqref{eq:coloring to partition}.
For $n\in\mathbb{N}$, we write $\mathop{[k]^{[n]}} $ to denote the
set of
$k$-colorings of $[n]$,
$\mathop{\mathcal{P}_{[n]}} $ to denote the set of partitions of
$[n]$, and
$\mathop{\mathcal{P}_{[n]:k}} $ to denote the subset of partitions of
$[n]$ with
at most $k$ blocks.

Any one-to-one mapping $\varphi\dvtx[m]\rightarrow[n]$, $m\leq n$,
determines a map $\mathop{[k]^{[n]}} \rightarrow[k]^{[m]}$, $x\mapsto
x^{\varphi
}$, where
%
%
\begin{equation}
\label{eq:composite mapping} x^{\varphi}=x^{\varphi(1)}\cdots x^{\varphi(m)}.
\end{equation}
We call the image in \eqref{eq:composite mapping} a \emph{composite
mapping} because $x\mapsto x^{\varphi}$ can be obtained by composing
the relabeling and restriction operations in \eqref{eq:relabeling} and
\eqref{eq:restriction}.
Let $\Rmn$ denote the restriction map $\mathop{[k]^{[n]}} \rightarrow
[k]^{[m]}$,
that is, $\Rmn x=x^{[m]}$.
To any one-to-one map $\varphi\dvtx[m]\rightarrow[n]$, there exists a
permutation $\sigma\dvtx[n]\rightarrow[n]$ such that $x^{\varphi
}=\Rmn
(x^{\sigma})$, relabeling by $\sigma$ followed by restriction to $[k]^{[m]}$.

For a partition $\pi\in\mathop{\mathcal{P}_{[n]}} $, relabeling,
restriction and
composite operations are defined by $\pi\mapsto\pi^{\sigma}$, $\pi
\mapsto\pi^{[m]}$, and $\pi\mapsto\pi^{\varphi}$, respectively, where
\begin{eqnarray*}
i&\sim&_{\hspace*{-2pt}\pi^{\sigma}}j\quad \Longleftrightarrow\quad\sigma(i)\sim_{\pi
}
\sigma(j),
\\
i&\sim&_{\hspace*{-2pt}\pi^{[m]}}j\quad \Longleftrightarrow\quad i\sim_{\pi} j\quad \mbox{and}
\\
i&\sim&_{\hspace*{-2pt}\pi^{\varphi}}j\quad \Longleftrightarrow\quad\varphi(i)\sim_{\pi
}
\varphi(j).
\end{eqnarray*}
When convenient, we abuse notation and also write $\Rmn$ to denote the
restriction $\mathop{\mathcal{P}_{[n]}} \rightarrow\mathcal
{P}_{[m]}$, that is, $\Rmn
\pi=\pi^{[m]}$, so that $\pi^{\varphi}=\Rmn(\pi^{\sigma})$ for some
$\sigma\dvtx[n]\rightarrow[n]$.

Any finite $k$-coloring can be embedded into a $k$-coloring of $\Nb$,
and likewise for partitions.
A $k$-coloring of $\Nb$ is an infinite $[k]$-valued sequence
$x=x^1x^2\cdots$ and is determined by its sequence of finite
restrictions $(x^{[1]},x^{[2]},\ldots)$.
A partition of $\Nb$ is defined similarly as a sequence of finite
partitions $(\pi^{[1]},\pi^{[2]},\ldots)$ for which $\pi^{[m]}=\Rmn
\pi
^{[n]}$, for every $m\leq n$.
As for finite sets, we denote $k$-colorings of $\Nb$ by $[k]^{\mathbb{N}}$,
partitions of $\Nb$ by $\partitionsN$, and partitions of $\Nb$ with at
most $k$ blocks by $\partitionsNk$.

For each $n\in\Nb$, $\Rn$ denotes the restriction map $[k]^{\mathbb{N}}
\rightarrow\mathop{[k]^{[n]}} $, or $\partitionsN\rightarrow\mathop
{\mathcal{P}_{[n]}} $.
The projective nature of both $[k]^{\mathbb{N}}$ and $\partitionsN$
endows each with a natural product-discrete topology.
With $\lambda,\lambda'$ denoting objects both in either $[k]^{\mathbb
{N}}$ or
$\partitionsN$, we define the ultrametric $d$ by
%
%
\begin{equation}
\label{eq:partition metric}d\bigl(\lambda,\lambda '\bigr):=2^{-n(\lambda,\lambda')},
\end{equation}
where $n(\lambda,\lambda'):=\max\{n\in\mathbb{N}\dvtx\Rn\lambda
=\Rn
\lambda
'\}$.
Under \eqref{eq:partition metric}, both $[k]^{\mathbb{N}}$ and
$\partitionsN$ are compact, separable and, therefore, Polish, metric spaces.
We equip $[k]^{\mathbb{N}}$ and $\partitionsNk$ with their discrete
$\sigma
$-fields, $\sigma \langle\bigcup_{n=1}^\infty\mathop
{[k]^{[n]}}  \rangle
$ and $\sigma \langle\bigcup_{n=1}^\infty\mathop{\mathcal
{P}_{[n]}}  \rangle
$, respectively.

\subsection{Exchangeability}\label{section:exch}
An infinite sequence $X:=(X_1,X_2,\ldots)$ of random variables is
called \emph{exchangeable} if its law is invariant under finite
permutations of its indices, that is, for each $n\in\mathbb{N}$,
\[
(X_{\sigma(1)},\ldots,X_{\sigma(n)})\equalinlaw(X_1,\ldots,X_n)\qquad\mbox{for every }\sigma\in\mathop{\mathscr{S}_{n}},
\]
where $\mathop{\mathscr{S}_{n}} $ denotes the symmetric group of
permutations of $[n]$.
By de Finetti's theorem (see, e.g., Aldous \cite
{AldousExchangeability}), the law of any exchangeable sequence $X\in
[k]^{\mathbb{N}}$ is determined by a unique directing probability
measure $\nu$ on the $(k-1)$-dimensional simplex
\[
\simplexk:= \Biggl\{(s_1,\ldots,s_k)\dvtx
s_i\geq0,\sum_{i=1}^k
s_i=1 \Biggr\}.
\]
In particular, conditional on $s\sim\nu$, $X_1,X_2,\ldots$ are
independent and identically distributed according to
\[
\mathbb{P}_s\{X_1=j\}=s_j,\qquad j=1,\ldots, k.
\]

A random partition $\Pi$ is exchangeable if $\Pi\equalinlaw\Pi
^{\sigma
}$ for all $\sigma\in\symmetricN$, where $\symmetricN$ is the set of
\emph{finite} permutations of $\Nb$, that is, permutations $\sigma
\dvtx\Nb
\rightarrow\Nb$ that fix all but finitely many elements.
Through \eqref{eq:coloring to partition}, any exchangeable $[k]$-valued
sequence $X$ projects to an exchangeable random partition $\Pi:=\Bb(X)$.
This construction of $\Pi$ is a special case of Kingman's paintbox
representation for exchangeable random partitions of $\mathbb{N}$
\cite
{Kingman1978b}.
If $X$ is directed by $\nu$, then we denote the law of $\Pi=\Bb(X)$
by~$\varrho_{\nu}$, the \emph{paintbox measure} directed by $\nu$.

With $f_i(X)$ defined in \eqref{eq:asymptotic frequency}, the
asymptotic frequency $|X|=(f_1(X),\ldots,\break  f_k(X))$ of any exchangeable
$k$-coloring exists almost surely.
Likewise for the asymptotic frequency of an exchangeable partition $\Pi
$, denoted $|\Pi|^{\downarrow}$, the vector of asymptotic block
frequencies listed in decreasing order of size which lives in the
ranked $k$-simplex \mbox{$\Delta_k^{\downarrow}:=\{(s_1,\ldots,s_k)\dvtx
s_1\geq
\cdots
\geq s_k\geq0,\mbox{ }\sum_is_i=1\}$}.

%
\begin{rmk}
To avoid measurability concerns, we can add the point $\partial$ to
both $\simplexk$ and $\Delta_k^{\downarrow}$ and put $|x|=\partial$
(resp., $|\pi|^{\downarrow}=\partial$) whenever the asymptotic
frequency of
$x\in[k]^{\mathbb{N}}$ (resp., $\pi\in\partitionsNk$) does not
exist. We
equip $\simplexk$, respectively, $\Delta_k^{\downarrow}$, with the
$\sigma
$-field generated by $|\cdot|\dvtx[k]^{\mathbb{N}}\rightarrow
\simplexk\cup\{
\partial\}$ and $|\cdot|^{\downarrow}\dvtx\partitionsNk\rightarrow
\Delta_k^{\downarrow}\cup\{\partial\}$, respectively. Beyond this point,
issues of measurability never arise, and so neither does the above formalism.
\end{rmk}

\subsection{Exchangeable Markov processes}

Let $\Xb=(X_t, t\in T)$ be a random collection in $[k]^{\mathbb{N}}$, with
$T$ either $\mathbb{Z}_+=\{0,1,\ldots\}$ (discrete-time) or $\mathbb
{R}_+=[0,\infty)$ (continuous-time).
We say $\Xb$ is \emph{Markovian} if, for every $t,t'\geq0$, the
conditional law of $X_{t+t'}$, given $\mathcal{F}_t:=\sigma\langle
X_s, s\leq t\rangle$, depends only on $X_t$ and $t'$.
Specifically, we distinguish between collections with finitely many
jumps in bounded intervals (Markov chains) and those with infinitely
many jumps in bounded intervals (Markov processes).
When speaking generally, we use the terminology and notation of Markov
processes as a catch-all.

The \emph{Markov semigroup} $\Pb=(\Pb_t, t\in T)$ of $\Xb=(X_t,
t\in
T)$ is defined for all bounded, measurable functions $g\dvtx
[k]^{\mathbb{N}}
\rightarrow\mathbb{R}$ by
%
%
\begin{equation}
\label{eq:semigroup} \Pb_tg(x):=\mathbb{E}_x
g(X_t),\qquad t\in T,
\end{equation}
the conditional expectation of $g(X_t)$ given $X_0=x$.
We say $\Xb$ \emph{enjoys the Feller property}, or is a \emph{Feller
process}, if for every bounded, continuous $g\dvtx[k]^{\mathbb
{N}}\rightarrow
\mathbb{R}$, its semigroup $\mathbf{P}$ satisfies:
\begin{itemize}
\item$\lim_{t\downarrow0}\mathbf{P}_tg(x)=g(x)$ for all $x\in
[k]^{\mathbb{N}}
$ and
\item$x\mapsto\mathbf{P}_tg(x)$ is continuous for all $t\in T$.
\end{itemize}

In general, since each $\Rn\dvtx[k]^{\mathbb{N}}\rightarrow\mathop
{[k]^{[n]}}
$ is a
many-to-one function, the restriction $\Xb^{[n]}$ need not be Markovian.
Under the product-discrete topology induced by~\eqref{eq:partition
metric}, exchangeability and consistency are equivalent to
exchangeability and the Feller property, and so we use the terms \emph
{consistency} and \emph{Feller} interchangeably.

%
\begin{prop}\label{prop:Feller equiv}
The following are equivalent for a Markov process $\bolds{\Lambda}$ on
either $[k]^{\mathbb{N}}$ or $\partitionsNk$:
\begin{longlist}[(ii)]
\item[(i)] $\bolds{\Lambda}$ is exchangeable and consistent under
subsampling.
\item[(ii)] $\bolds{\Lambda}$ is exchangeable and enjoys the Feller
property.
\end{longlist}
\end{prop}

\subsection{Coset decompositions and associated mappings}\label
{section:cosets}
For fixed $k\in\Nb$, we define the \emph{coset decomposition} of
$x\in
[k]^{\mathbb{N}}$ by the $k$-tuple $(x_1,\ldots,x_k)$, where
%
%
\begin{equation}
\label{eq:coset} x_i=x^{i}x^{i+k}x^{i+2k}
\cdots, \qquad i=1,\ldots,k.
\end{equation}
In words, the \emph{$i$th coset} of $x$ is the subsequence of $x$
including every $k$th element, beginning at coordinate $i$.
Through \eqref{eq:coset}, the sets $[k]^{\mathbb{N}}$ and
$[k]^{\mathbb{N}\otimes k}\cong
[k]^{\mathbb{N}}\times\cdots\times[k]^{\mathbb{N}}$ ($k$ times)
are in one-to-one
correspondence, but we sometimes prefer one representation over the other.
To distinguish between representations, we write:
\begin{itemize}
\item$x=x^1x^2\cdots$ to denote the object in $[k]^{\mathbb{N}}$ and
\item$x=(x_1,\ldots,x_k)$ to denote the coset representation in
$[k]^{\mathbb{N}\otimes k}$, with each coset written
\[
x_i=x_i^1x_i^2
\cdots=x^ix^{i+k}\cdots.
\]
\end{itemize}
We usually write $x$ to denote an object \emph{initially} defined in
$[k]^{\mathbb{N}}$ and $M$ to denote an object \emph{initially}
defined in
$[k]^{\mathbb{N}\otimes k}$.
The importance of this decomposition becomes apparent in Section~\ref
{section:random matrix product}.

For $n\in\Nb$, the restriction of $M\in[k]^{\mathbb{N}\otimes k}$
to $\mathop
{[k]^{[n]\otimes k}} \cong
\mathop{[k]^{[n]}} \times\cdots\times\mathop{[k]^{[n]}} $ ($k$~times) is defined
componentwise by
%
%
\begin{equation}
\label{eq:coset restriction} M^{[n]}:=\bigl(M_1^{[n]},
\ldots,M_k^{[n]}\bigr).
\end{equation}
Likewise, a $k$-tuple of finite permutations $\sigma_1,\ldots,\sigma
_k\dvtx\Nb\rightarrow\Nb$ acts on $M\in[k]^{\mathbb{N}}$ by
%
%
\begin{equation}
\label{eq:coset relabeling} M^{\sigma_1,\ldots,\sigma_k}:=\bigl(M_1^{\sigma_1},
\ldots,M_k^{\sigma_k}\bigr).
\end{equation}

Any $M\in\mathop{[k]^{[n]\otimes k}} $ functions as a map $\mathop
{[k]^{[n]}} \rightarrow\mathop{[k]^{[n]}} $.
For each $x\in\mathop{[k]^{[n]}} $, we define the injection $\varphi
_x\dvtx[n]\rightarrow[nk]$ by
%
%
\begin{equation}
\label{eq:varphi-x} \varphi_x(j):=x^j+(j-1)k,\qquad j=1,\ldots,n.
\end{equation}
For any $M\in[k]^{\mathbb{N}\otimes k}$, its restriction $M^{[n]}$ to
$\mathop
{[k]^{[n]\otimes k}} $, as
in \eqref{eq:coset restriction}, is in correspondence with a unique
$k$-coloring $M^1\cdots M^{nk}$ of $[nk]$.
Using \eqref{eq:composite mapping}, we define $M^{[n]}\dvtx\mathop{[k]^{[n]}}
\rightarrow\mathop{[k]^{[n]}} $ by
%
%
\begin{equation}
\label{eq:coset mapping} M^{[n]}(x):=M^{\varphi_x}=M^{x^1}M^{x^2+k}
\cdots M^{x^n+(n-1)k},\qquad x\in\mathop{[k]^{[n]}}.
\end{equation}
The finite maps $(M^{[n]}, n\in\Nb)$ derived from $M$ determine a
unique map $M\dvtx[k]^{\mathbb{N}}\rightarrow[k]^{\mathbb{N}}$.

Importantly, each $M\in[k]^{\mathbb{N}\otimes k}$ determines a
Lipschitz continuous map
in the metric \eqref{eq:partition metric}.
The identity map $\idk\dvtx[k]^{\mathbb{N}}\rightarrow[k]^{\mathbb
{N}}$ corresponds to
the infinite repeating pattern $12\cdots k$,
%
%
\begin{equation}
\label{eq:canonk} \canonk=12\cdots k12\cdots k\cdots,
\end{equation}
for example, $\mathbf{Z}_2=121212\cdots,$ $\mathbf{Z}_3=123123\cdots,$ and
so on.
The coset decomposition of $\canonk$ is $(\mathbf{1},\mathbf
{2},\ldots,\mathbf{k})$, where $\mathbf{i}=\mathit{iii}\cdots$ is the infinite sequence of
all $i$'s, for each $i=1,2,\ldots.$
For $n\in\Nb$, we write $\canonkn$ to denote the restriction of
$\canonk
$ to $\mathop{[k]^{[n]\otimes k}} $ and $\idnk$ to denote its
associated identity map
$\mathop{[k]^{[n]}} \rightarrow\mathop{[k]^{[n]}} $.
By definition \eqref{eq:coset relabeling}, $\canonk$, and hence
$\canonkn$, is invariant under relabeling by any $k$-tuple $\sigma
_1,\ldots,\sigma_k$ of permutations.

Since any mapping $M\dvtx[k]^{\mathbb{N}}\rightarrow[k]^{\mathbb
{N}}$ is determined by
its coset decomposition $(M_1,\ldots,M_k)$, we can define the \emph
{asymptotic frequency} of $M$ by the $k$-tuple $(|M_1|,\ldots,|M_k|)$,
provided each $|M_i|$ exists.
We express the asymptotic frequency of $M$ as a stochastic matrix
$|M|_k=S=(S_{\mathit{ii}'}, 1\leq i,i'\leq k)$, where
%
%
\begin{equation}
\label{eq:asymptotic frequency matrix} S_{\mathit{ii}'}:=\lim_{n\rightarrow\infty}n^{-1}\sum
_{j=1}^{n}\mathbf{1}\bigl\{
M_i^{j}=i'\bigr\},\qquad 1\leq i,i'
\leq k.
\end{equation}

\section{Discrete-time cut-and-paste chains}\label{section:random
matrix product}
In this section, $\Xb=(X_m, m\geq0)$ denotes a discrete-time
exchangeable and consistent Markov chain on $[k]^{\mathbb{N}}$, and
$\Xb
^{[n]}=(X_m^{[n]}, m\geq0)$ its restriction to $\mathop{[k]^{[n]}} $,
for each
$n=1,2,\ldots.$
By assumptions~(A) and (B), each $\Xb^{[n]}$ is an exchangeable Markov
chain with transition probability measure
\[
P_n\bigl(x,x'\bigr):=\mathbb{P}\bigl
\{X_1=x'\mid X_0=x\bigr\},\qquad
x,x'\in\mathop{[k]^{[n]}}.
\]
Exchangeability implies $P_n(x,x')=P_n(x^{\sigma},x'^{\sigma})$ for all
permutations $\sigma\dvtx[n]\rightarrow[n]$, while consistency relates
$(P_n, n\in\Nb)$ through
\[
P_m\bigl(x,x'\bigr)=P_n\bigl(x^*,
\mathbf{R}^{-1}_{m,n}\bigl(x'\bigr)\bigr),\qquad
x,x'\in[k]^{[m]},
\]
for all $x^*\in\mathbf{R}^{-1}_{m,n}(x)=\{\hat{x}\in\mathop
{[k]^{[n]}} \dvtx\hat
{x}^{[m]}=x\}$.
Writing $P$ to denote the transition probability measure of $\Xb$ on
$[k]^{\mathbb{N}}$, we conclude
%
%
\begin{equation}
\label{eq:induced fidi tps} P_n\bigl(x,x'\bigr)=P\bigl(x^*,
\mathbf{R}^{-1}_n\bigl(x'\bigr)\bigr),\qquad
x,x'\in\mathop{[k]^{[n]}}, \mbox{for all }x^*\in
\mathbf{R}^{-1}_n(x),
\end{equation}
for every $n\in\Nb$.

Theorem~\ref{thmm:cut-and-paste chain} asserts that $P$ is determined
by a unique probability measure $\Sigma$ on $\stochk$.
We construct $\Sigma$ directly from $P$ using the connection between
$k$-colorings and stochastic matrices from Section~\ref{section:cosets}.
For $\canonk$ in \eqref{eq:canonk}, we define a probability measure
$\chi$ on $[k]^{\mathbb{N}\otimes k}$ by
%
%
\begin{equation}
\label{eq:chi} \chi(\cdot):=P(\canonk,\cdot).
\end{equation}

%
\begin{defn}[(Coset exchangeability)]\label{defn:coset exchangeability}
A random mapping $M=(M_1,\break \ldots, M_k)\in[k]^{\mathbb{N}\otimes k}$ is
\emph{coset
exchangeable} if
%
%
\begin{equation}
\label{eq:coset exchangeable} (M_1,\ldots,M_k)\equalinlaw
\bigl(M_1^{\sigma_1},\ldots,M_k^{\sigma_k}\bigr)\qquad
\mbox{for all }\sigma_1,\ldots,\sigma_k\in\symmetricN.
\end{equation}
\end{defn}

For any random mapping $M$ constructed from a random $k$-coloring\break
through~\eqref{eq:coset}, exchangeability implies coset
exchangeability, but not the reverse.
By assumption, $P$ is an exchangeable transition probability on
$[k]^{\mathbb{N}}$ and the coset decomposition of $\canonk$ is invariant
under coset relabeling \eqref{eq:coset relabeling}; hence, $\chi$
defined in~\eqref{eq:chi} is coset exchangeable and the asymptotic
frequency of $M\sim\chi$, as defined in \eqref{eq:asymptotic frequency
matrix}, exists with probability one.
We denote the law of $|M|_k$ by $|\chi|_k$.

We complete the proof of Theorem~\ref{thmm:cut-and-paste chain} by
showing that a random $k$-coloring $X'$ generated by first drawing
$M\sim\chi$ and then putting $X'=M(x)$, for fixed $x\in[k]^{\mathbb
{N}}$, is
a draw from $P(x,\cdot)$.
By consistency, we need only show that $M^{[n]}(x)\sim P_n(x,\cdot)$
for every $x\in\mathop{[k]^{[n]}} $, for every $n\in\Nb$.
We have defined $\canonk$ so that
\[
\canonk(x)=\mathbf{Z}_{k,n}^{\varphi_x}=\mathbf{Z}_k^{x^1}
\cdots \mathbf {Z}_{k}^{x^n+(n-1)k}=x^1\cdots
x^n=x \qquad\mbox{for all }x\in\mathop {[k]^{[n]}}.
\]
By \eqref{eq:induced fidi tps} and \eqref{eq:chi}, the restriction of
$M\sim\chi$ to $\mathop{[k]^{[n]\otimes k}} $ is distributed as
\[
M^{[n]}\sim\chi^{(n)}(\cdot)=P_{nk}(\canonkn,\cdot),
\]
which combines with \eqref{eq:coset mapping} to imply $M^{[n]}(x)\sim
P_n(x,\cdot)$.

We have proven the following prelude to Theorem~\ref
{thmm:cut-and-paste chain}.

%
\begin{thmm}\label{thmm:characteristic measure discrete}
Let $\Xb=(X_m, m\geq0)$ be a discrete-time, exchangeable, consistent
Markov chain on $[k]^{\mathbb{N}}$.
Then there exists a probability measure $\chi$ on $[k]^{\mathbb
{N}\otimes k}$ such that
$\Xb^*=(X_m^*, m\geq0)$ is a version of $\Xb$, where
$X^*_0\equalinlaw
X_0$ and
\[
X_m^*=(M_m\circ\cdots\circ M_1)
\bigl(X_0^*\bigr),\qquad m\geq1,
\]
for $M_1,M_2,\ldots$ drawn i.i.d. from $\chi$.
\end{thmm}

To establish Theorem~\ref{thmm:cut-and-paste chain}, we must show that
$\chi$ is determined by a unique probability measure on $\stochk$.
By \eqref{eq:asymptotic frequency matrix} and coset exchangeability,
$\chi$ induces a probability measure $|\chi|_k$ on $\stochk$.
By de Finetti's theorem, the components of $M\sim\chi$, given
$|M|_k=S$, are conditionally independent with distribution
%
%
\begin{equation}
\label{eq:independent components} \Pbb_S\bigl\{M^{i+(j-1)k}=i'\bigr
\}=S_{\mathit{ii}'}, \qquad i,i'=1,\ldots,k; j=1,2,\ldots.
\end{equation}
We write $\mu_S$ to denote the conditional distribution of $M$, given
$|M|_k=S$, as in~\eqref{eq:independent components} and
%
%
\begin{equation}
\label{eq:mu-mixture} \mu_{\Sigma}(\cdot):=\int_{\stochk}
\mu_S(\cdot)\Sigma(dS)
\end{equation}
to denote the mixture of $\mu_S$-measures with respect to $\Sigma$.
By \eqref{eq:independent components}, the components $Y^1Y^2\cdots$ of
$M(x)$ are conditionally independent given $|M|_k=S$ and have distribution
%
%
\begin{equation}
\label{eq:conditional tps} \Pbb_S\bigl\{Y^j=i'\mid
x^j=i\bigr\}=S_{\mathit{ii}'},\qquad j=1,2,\ldots.
\end{equation}
For every $n\in\Nb$, the unconditional law of $M^{[n]}(x)$ is thus
\[
P_n\bigl(x,x'\bigr)=\int_{\stochk}
\prod_{j=1}^n S\bigl(x^j,x'^j
\bigr)|\chi|_k(dS),\qquad x'\in\mathop{[k]^{[n]}}.
\]
Putting $\Sigma:=|\chi|_k$ establishes Theorem~\ref
{thmm:cut-and-paste chain}.

%
\begin{rmk}
We call $\Xb^*_{\Sigma}$ in Theorem~\ref{thmm:cut-and-paste chain} an
\emph{(exchangeable) cut-and-paste chain} with directing measure
$\Sigma
$ and cut-and-paste measure $\mu_{\Sigma}$.
\end{rmk}

From Theorems \ref{thmm:cut-and-paste chain} and \ref
{thmm:characteristic measure discrete}, we can generate a version of
$\Xb$ by drawing $X_0$ from the initial distribution of $\Xb$ and
$M_1,M_2,\ldots$ i.i.d. from $\mu_{\Sigma}$.
Given $X_0,M_1,M_2,\ldots,$ we define
%
%
\begin{equation}
\label{eq:iterated action} X_m:=M_m(X_{m-1})=(M_m
\circ\cdots\circ M_1) (X_0),\qquad m\geq1.
\end{equation}
By de Finetti's theorem, $|X_0|=(f_1(X_0),\ldots,f_k(X_0))$ exists
almost surely and $|M_1|_k,|M_2|_k,\ldots$ is an i.i.d. sequence from
$\Sigma$.
By the construction of $\Xb$ in \eqref{eq:iterated action}, $X_1$ is
chosen from the conditional transition probability in \eqref
{eq:conditional tps}, with $S=|M_1|_k$.
By the strong law of large numbers, $f_{i'}(X_1)$ exists almost surely
for every $i'=1,\ldots,k$ and equals the $i'$th component of
$|X_0|S_1$, that is, %
\[
f_{i'}(X_1)=\sum_{i=1}^kf_i(X_0)S_1
\bigl(i,i'\bigr).
\]
By induction, the components of $|X_m|$, given $|X_{m-1}|$ and
$|M_m|_k$, equal
\[
|X_{m-1}||M_m|_k=|X_0||M_1|_k
\cdots|M_m|_k \qquad\mbox{for every }m\geq1,
\]
and Theorem~\ref{thmm:induced simplex chain} follows.

\section{Continuous-time cut-and-paste processes}\label{section:Levy-Ito}

We now let $\Xb=(X_t, t\geq0)$ denote an exchangeable, consistent
Markov process in continuous-time.
We have noted previously that $\Xb$ can jump infinitely often in
bounded intervals, but its finite restrictions can jump only finitely often.
To characterize the behavior of $\Xb$, we use a Poisson point process
to build a version sequentially through its finite restrictions.
Similar to our discrete-time construction \eqref{eq:iterated action},
we define the intensity measure of the Poisson point process directly
from the transition law of $\Xb$.
Dissimilar to the discrete-time case, this intensity need not be finite.

Let $\chi$ be a coset exchangeable measure on $[k]^{\mathbb{N}\otimes
k}$ satisfying
%
%
\begin{eqnarray}
\label{eq:regularity chi} \chi\bigl(\{\idk\}\bigr)=0 \quad\mbox{and}\quad \chi\bigl(\bigl\{M
\in[k]^{\mathbb{N}\otimes
k}\dvtx M^{[n]}\neq \idnk\bigr\}\bigr)<\infty
\nonumber
\\[-8pt]
\\[-8pt]
\eqntext{\mbox{for
all }n\in\Nb.}
\end{eqnarray}
We construct a process $\Xb^*_{\chi}=(X_t^*, t\geq0)$ through its
finite restrictions $(\Xb^{*[n]}_{\chi},\break  n\in\Nb)$ as follows.
Let $\Mb=\{(t,M_t)\}\subseteq\mathbb{R}_+\times[k]^{\mathbb
{N}\otimes k}$ be a Poisson
point process with intensity $dt\otimes\chi$, where $dt$ denotes
Lebesgue measure on $[0,\infty)$.
Given an exchangeable initial state $X_0\in[k]^{\mathbb{N}}$, we put
$X_0^{*[n]}=X_0^{[n]}$ and, for each $t>0$:
\begin{itemize}
\item if $t>0$ is an atom time of $\Mb$ for which $M_t^{[n]}\neq\idnk$,
we put $X^{*[n]}_t:=M_{t}^{[n]}(X_{t-}^{*[n]})$,
\item otherwise, we put $X_t^{*[n]}=X_{t-}^{*[n]}$.
\end{itemize}

This construction of each $\Xb^{*[n]}_{\chi}$ is a continuous-time
analog to the discrete-time construction in \eqref{eq:iterated action};
it differs only in the random time between jumps and the possibility of
infinitely many jumps in the limiting process.
We have constructed each $\Xb^{*[n]}_{\chi}$ from the same Poisson
process so that $(\Xb^{*[n]}_{\chi}, n\in\Nb)$ is compatible, that
is, $X_t^{*[m]}=\Rmn X_t^{*[n]}$ for all $t\geq0$ and $m\leq n$, and
determines a unique $[k]^{\mathbb{N}}$-valued process $\Xb^*_{\chi}$.

%
\begin{prop}
Let $\chi$ be a coset exchangeable measure on $[k]^{\mathbb{N}\otimes
k}$ that
satisfies \eqref{eq:regularity chi}, and let $\Xb^*_{\chi}$ be as
constructed from the Poisson point process $\Mb$ with intensity
$dt\otimes\chi$.
Then $\Xb^*_{\chi}$ is an exchangeable, consistent Markov process on
$[k]^{\mathbb{N}}$.
\end{prop}

\begin{pf}
For each $n\in\Nb$, $\mathbf{X}^{*{[n]}}_{\chi}$ is a Markov chain by
assumption \eqref{eq:regularity chi} and its Poisson point process
construction.
Moreover, $\Rmn X^{*[n]}_t=X^{*[m]}_t$ for all $t\geq0$, for all
$m\leq
n$, and so $(\Xb^{*[n]}_{\chi}, n\in\Nb)$ determines a unique Markov
process $\Xb^{*}_{\chi}$ on $[k]^{\mathbb{N}}$.
Exchangeability of $\Xb^*_{\chi}$ follows by coset exchangeability of
$\chi$, since all of its finite restrictions to $\mathop
{[k]^{[n]\otimes k}} $ are finite,
coset exchangeable measures.
\end{pf}
%
%
\begin{cor}
Every coset exchangeable measure $\chi$ on $[k]^{\mathbb{N}\otimes
k}$ satisfying~\eqref
{eq:regularity chi} determines the jump rates of an exchangeable Feller
process on~$[k]^{\mathbb{N}}$.
\end{cor}

A measure satisfying \eqref{eq:regularity chi} can be constructed
directly from the transition rates of $\Xb$.
By assumption, each finite restriction $\Xb^{[n]}=(X_t^{[n]}, t\geq0)$
is a c\`adl\`ag, exchangeable Markov\vadjust{\goodbreak} process on $\mathop{[k]^{[n]}} $.
Since $\mathop{[k]^{[n]}} $ is finite, the evolution of $\Xb^{[n]}$ is
characterized by its jump rates
%
%
\begin{equation}
\label{eq:finite jump rates} Q_n\bigl(x,x'\bigr):=\lim
_{t\downarrow0}\frac{1}{t}\mathbb {P}\bigl(X^{[n]}_t=x'
\mid X^{[n]}_0=x\bigr),\qquad x\neq x'\in
\mathop{[k]^{[n]}},
\end{equation}
which satisfy
%
%
\begin{equation}
\label{eq:finite markov rates} Q_n\bigl(x,\mathop{[k]^{[n]}} \setminus\{x\}
\bigr)<\infty\qquad\mbox{for all }x\in \mathop{[k]^{[n]}},
\end{equation}
are exchangeable in the sense that, for every $\sigma\in\mathop
{\mathscr{S}_{n}} $,
%
%
\begin{equation}
\label{eq:exch markov rates} Q_n\bigl(x,x'\bigr)=Q_n
\bigl(x^{\sigma},x'^{\sigma}\bigr),\qquad x\neq
x'\in\mathop {[k]^{[n]}},
\end{equation}
and are consistent,
%
%
\begin{eqnarray}
\label{eq:consistent markov rates} Q_m\bigl(x,x'\bigr)=Q_n
\bigl(x^*,\mathbf{R}_{m,n}^{-1}\bigl(x'\bigr)
\bigr),
\nonumber
\\[-8pt]
\\[-8pt]
\eqntext{ x\neq x'\in [k]^{[m]}, \mbox{for all } x^*\in
\mathbf{R}^{-1}_{m,n}(x).}
\end{eqnarray}
For each $n\in\Nb$, we define
%
%
\begin{equation}
\label{eq:finite chi cont} \chi_n(M):=Q_n(\canonkn,M),\qquad M\in
\mathop{[k]^{[n]\otimes k}} \setminus\{\idnk\}.
\end{equation}

%
\begin{lemma}\label{lemma:existence matrix measure}
The collection $(\chi_n, n\in\mathbb{N})$ in \eqref{eq:finite chi
cont} is coset exchangeable and satisfies
\[
\chi_m(M)=\chi_n\bigl(\bigl\{M^*\in
\mathop{[k]^{[n]\otimes k}} \dvtx M^{*[m]}=M\bigr\} \bigr)\qquad\mbox{for all
} M\in[k]^{[m]\otimes k},
\]
for all $m\leq n$.
\end{lemma}
\begin{pf}
This follows from the definition of $\chi_n$ in \eqref{eq:finite chi
cont}, the correspondence $[k]^{\mathbb{N}}\leftrightarrow
[k]^{\mathbb{N}\otimes k}$ in
\eqref
{eq:coset}, and conditions \eqref{eq:finite markov rates}, \eqref
{eq:exch markov rates} and \eqref{eq:consistent markov rates}.
\end{pf}

%
\begin{prop}\label{prop:caratheodory}
Let $(\chi_n, n\in\mathbb{N})$ be defined in \eqref{eq:finite chi cont}.
Then there exists a unique coset exchangeable measure $\chi$ on
$[k]^{\mathbb{N}\otimes k}$ satisfying \eqref{eq:regularity chi} and
\begin{eqnarray}
\chi\bigl(\bigl\{M^*\in[k]^{\mathbb{N}\otimes k}\dvtx M^{*[n]}=M\bigr\}\bigr)=
\chi _n(M),\nonumber\\
 \eqntext{M\in\mathop {[k]^{[n]\otimes k}} \setminus\{\idnk\},
\mbox{ for every }n\in\mathbb{N}.}
\end{eqnarray}
\end{prop}
\begin{pf}
Because $\bigcup_{n=1}^\infty\mathop{[k]^{[n]\otimes k}} $ is a
generating $\pi$-system of
the product $\sigma$-field over $[k]^{\mathbb{N}\otimes k}$, we need
only determine
$\chi$ on subsets of the form
\[
\bigl\{M^*\in[k]^{\mathbb{N}\otimes k}\dvtx M^{*[n]}=M\bigr\},
\]
for every $n\in\mathbb{N}$ and $M\in\mathop{[k]^{[n]\otimes k}} $.
Lemma~\ref{lemma:existence matrix measure} implies
\[
\chi_m(M)=\chi_n\bigl(\bigl\{M^*\in
\mathop{[k]^{[n]\otimes k}} \dvtx M^{*[m]}=M\bigr\} \bigr)=\sum
_{M^*\in\mathop{[k]^{[n]\otimes k}}
\dvtx M^{*[m]}=M}\chi_n\bigl(M^*\bigr),
\]
for all $m\leq n$ and $M\in[k]^{[m]\otimes k}$.
Therefore, $\chi$ defined by
%
%
\begin{equation}
\label{eq:finite chi}\qquad \chi\bigl(\bigl\{M^*\in[k]^{\mathbb{N}\otimes k}\dvtx M^{*[n]}=M
\bigr\}\bigr)=\chi _n(M), \qquad M\in [k]^{[n]\otimes k}\setminus\{
\operatorname{id}_{k,n}\},
\end{equation}
is additive, and Caratheodory's extension theorem implies $\chi$ has a
unique extension to a measure on $[k]^{\mathbb{N}\otimes k}\setminus\{
\idk\}$.

To satisfy the first half of \eqref{eq:regularity chi}, we simply put
$\chi(\{\idk\})=0$.
For the second half, \eqref{eq:finite markov rates} implies
\begin{eqnarray*}
\chi\bigl(\bigl\{M\in[k]^{\mathbb{N}\otimes k}\dvtx M^{[n]}\neq\idnk\bigr\}
\bigr)&=&\chi_n\bigl(\mathop {[k]^{[n]\otimes k}} \setminus\{ \idnk\}
\bigr)\\
&=&Q_{nk}\bigl(\canonkn,[k]^{[nk]}\setminus\{\canonkn\}
\bigr)<\infty.
\end{eqnarray*}
This completes the proof.
\end{pf}

The measure $\chi$ in Proposition~\ref{prop:caratheodory} ties the
Poissonian construction of $\Xb^*_{\chi}$ to $\Xb$, as the next
theorem shows.

%
\begin{thmm}\label{thmm:characteristic measure Feller}
Let $\Xb$ be a continuous-time, exchangeable, consistent Markov process
on $[k]^{\mathbb{N}}$.
Then there exists a coset exchangeable measure $\chi$ on $[k]^{\mathbb
{N}\otimes k}$
satisfying \eqref{eq:regularity chi} such that $\Xb_{\chi}^*$ is a
version of $\Xb$.
\end{thmm}

\begin{pf}
Let $\chi$ be the coset exchangeable measure with finite-dimensional
distributions \eqref{eq:finite chi cont}.
By Proposition~\ref{prop:caratheodory}, $\chi$ satisfies \eqref
{eq:regularity chi}.

Let $\mathbf{X}_{\chi}^*$ be the Markov process constructed from $\Mb$
with intensity $dt\otimes\chi$.
The total intensity at which events occur in $\mathbf{M}$ is $\chi
([k]^{\mathbb{N}\otimes k})$.
For $n\in\mathbb{N}$, the atom times of $\mathbf{X}^{*[n]}_{\chi}$ are
a thinned version of the atom times of $\mathbf{M}$.
In the construction of $\mathbf{X}^{*{[n]}}_{\chi}$, an atom
$(t,M_t)\in
\mathbf{M}$ results in a jump in $\mathbf{X}^{*{[n]}}_{\chi}$ if and
only if $M_{t}^{[n]}\neq\idnk$ and $M_{t}^{[n]}(X^{*[n]}_{t-})\neq
X^{*[n]}_{t-}$.
By the thinning property of Poisson processes, given
$X^{*[n]}_{t-}=x\in
\mathop{[k]^{[n]}} $, the total intensity at which $\Xb_{\chi
}^{*[n]}$ jumps
from state $x$ to $x'\neq x$ is ${\chi_n(\{M\in\mathop
{[k]^{[n]\otimes k}} \dvtx M(x)=x'\})}$.
And by \eqref{eq:exch markov rates} and \eqref{eq:consistent markov rates},
\begin{eqnarray*}
\chi_n\bigl(\bigl\{M\in\mathop{[k]^{[n]\otimes k}} \dvtx
M(x)=x'\bigr\}\bigr)&=&\sum_{M:M(x)=x'}Q_{nk}(
\canonkn,M)
\\
&=&Q_{nk}\bigl(\canonkn,\bigl\{z\in[k]^{[nk]}\dvtx
z^{\varphi_x}=x'\bigr\}\bigr)
\\
&=&Q_n\bigl(x,x'\bigr).
\end{eqnarray*}
It follows that the total intensity of jumps out of $x$ is
\[
\chi_n\bigl(\bigl\{M\in\mathop{[k]^{[n]\otimes k}} \dvtx M(x)\neq x
\bigr\} \bigr)=Q_n\bigl(x,\mathop{[k]^{[n]}} \setminus\{x\}
\bigr)<\infty,
\]
and, for each $n\in\mathbb{N}$, $\mathbf{X}^{*[n]}_{\chi}$ is an
exchangeable Markov process with jump rates $Q_n(\cdot,\cdot)$.
Kolmogorov's extension theorem implies $\mathbf{X}^*_{\chi}$ is a
version of $\mathbf{X}$.
\end{pf}

\subsection{L\'evy--It\^o representation}\label{subsec:Levy-Ito}

Our entire discussion climaxes in Theorem~\ref{thmm:Levy-Ito}, the L\'
evy--It\^o representation.
For any exchangeable, consistent Markov process on $[k]^{\mathbb{N}}$, its
characteristic measure $\chi$ has two unique components:
a measure $\Sigma$ on $k\times k$ stochastic matrices for which
%
%
\begin{equation}
\label{eq:regularity Sigma2} \Sigma\bigl(\{I_k\}\bigr)=0 \quad\mbox{and}\quad \int
_{\stochk}(1-S_*)\Sigma (dS)<\infty,
\end{equation}
where $S_*:=\min(S_{11},\ldots,S_{kk})$, and a collection $\mathbf
{c}=(\mathbf{c}_{\mathit{ii}'}, 1\leq i\neq i'\leq k)$ of nonnegative constants.

For $1\leq i\neq i'\leq k$ and $n\in\Nb$, we define $\rho^{(n)}_{\mathit{ii}'}$
as the point mass at $\kappa_{\mathit{ii}'}^{(n)}=(z_1,\ldots,z_k)\in
[k]^{\mathbb{N}\otimes k}$, where
\[
z_j^{j'}=\cases{ %
i',& \quad $j=i, j'=n,$
\vspace*{2pt}\cr
j,&\quad $\mbox{otherwise.}$ }
\]
In words, $\rho^{(n)}_{\mathit{ii}'}$ charges only the map $\kappa_{\mathit{ii}'}^{(n)}$
that fixes all but the $n$th coordinate of every $x\in[k]^{\mathbb
{N}}$: if
$x^n=i$, then the $n$th coordinate of $\kappa_{\mathit{ii}'}^{(n)}(x)$ is $i'$;
otherwise, the $n$th coordinate is also unchanged.
We call each $\kappa_{\mathit{ii}'}^{(n)}$ a \emph{single-index flip}.
For example, with $k=3$, $\rho^{(3)}_{12}$ puts unit mass at $\kappa
_{12}^{(3)}=(1121\cdots,2222\cdots,3333\cdots)$.
The measure
\[
\rho_{\mathit{ii}'}(\cdot):=\sum_{n=1}^{\infty}
\rho_{\mathit{ii}'}^{(n)}(\cdot),\qquad 1\leq i\neq i'\leq k,
\]
puts unit mass at every single-index flip from $i$ to $i'$.

For any $\Sigma$ satisfying \eqref{eq:regularity Sigma2} and any
collection $(\mathbf{c}_{\mathit{ii}'}, 1\leq i\neq i'\leq k)$ of nonnegative
constants, we define
%
%
\begin{equation}
\label{eq:chi-continuous} \chi_{\Sigma,\mathbf{c}}:=\mu_{\Sigma}+\sum
_{1\leq i\neq i'\leq
k}\mathbf{c}_{\mathit{ii}'}\rho_{\mathit{ii}'},
\end{equation}
where $\mu_{\Sigma}$ was defined in \eqref{eq:mu-mixture}.

%
\begin{prop}
Let $\Sigma$ satisfy \eqref{eq:regularity Sigma2} and $\mathbf
{c}=(\mathbf{c}_{\mathit{ii}'}, 1\leq i\neq i'\leq k)$ be nonnegative constants.
Then $\chi_{\Sigma,\mathbf{c}}$ defined in \eqref
{eq:chi-continuous} is
a coset exchangeable measure satisfying \eqref{eq:regularity chi}.
\end{prop}

\begin{pf}
We treat each term of $\chi_{\Sigma,\mathbf{c}}$ separately.

Clearly, $\mu_\Sigma(\{\idk\})=0$ by the first half of \eqref
{eq:regularity Sigma2} and the strong law of large numbers. Now, for
every $n\in\mathbb{N}$ and $S\in\stochk$, we have
\[
\mu_S\bigl(\bigl\{M\dvtx M^{[n]}\neq\idnk\bigr\}\bigr)\leq
\sum_{j=1}^k\mu_S\bigl(\bigl\{
M\dvtx M_{j}^{[n]}\neq \mathbf{j}^{[n]}\bigr\}\bigr)
\leq k\bigl(1-S_*^n\bigr)\leq nk (1-S_*),
\]
where $\mathbf{j}=jj\cdots\in[k]^{\mathbb{N}}$ and $\mathbf
{j}^{[n]}:=j\cdots
j$ is its restriction to $\mathop{[k]^{[n]}} $.
By \eqref{eq:regularity Sigma2},
\[
\mu_\Sigma\bigl(\bigl\{M\dvtx M^{[n]}\neq\idnk\bigr\}\bigr)
\leq nk\int_{\stochk
}(1-S_*)\Sigma (dS)<\infty.
\]

The first half of \eqref{eq:regularity chi} is satisfied by $\sum_{i\neq i'}\mathbf{c}_{\mathit{ii}'}\rho_{\mathit{ii}'}$ because each $\rho_{\mathit{ii}'}$
charges only single-index flips.
Furthermore, with $c^*:=\max_{1\leq i\neq i'\leq k}\mathbf
{c}_{\mathit{ii}'}<\infty$,
\begin{eqnarray*}
\sum_{1\leq i\neq i'\leq k}\mathbf{c}_{\mathit{ii}'}
\rho_{\mathit{ii}'}\bigl(\bigl\{M\dvtx M^{[n]}\neq \idnk\bigr\}\bigr)&\leq&
c^*\sum_{1\leq i\neq i'\leq k}\sum_{j=1}^{n}
\rho _{\mathit{ii}'}^{(j)}\bigl([k]^{\mathbb{N}\otimes k}\bigr)\\
&=&nk(k-1)c^*<
\infty.
\end{eqnarray*}
Thus, $\chi_{\Sigma,\mathbf{c}}$ satisfies \eqref{eq:regularity chi}.

Coset exchangeability of $\chi_{\Sigma,\mathbf{c}}$ follows since it is
the sum of coset exchangeable measures.
\end{pf}

Now, the denouement.

\begin{pf*}{Proof of Theorem~\ref{thmm:Levy-Ito}}
By Theorem~\ref{thmm:characteristic measure Feller}, every exchangeable
Feller process on $[k]^{\mathbb{N}}$ admits a version $\Xb^*_{\chi
}$, for
$\chi$ satisfying \eqref{eq:regularity chi}.
In Theorem~\ref{thmm:Levy-Ito}, we assert that $\chi$ can be decomposed
as in \eqref{eq:chi-continuous}.
To prove this, we proceed in three steps:
\begin{longlist}[(iii)]
\item[(i)] $\chi$-almost every $M\in[k]^{\mathbb{N}\otimes k}$
possesses asymptotic
frequency \mbox{$|M|_k\in\stochk$},
\item[(ii)] there exists a unique measure $\Sigma$ satisfying \eqref
{eq:regularity Sigma2} such that the restriction of $\chi$ to $\{M\in
[k]^{\mathbb{N}\otimes k}\dvtx |M|_k\neq I_k\}$ is a cut-and-paste measure,
\[
\mathbf{1}_{\{|M|_k\neq I_k\}}\chi(dM)=\mu_{\Sigma}(dM)\quad \mbox{and}
\]
\item[(iii)] there exist unique nonnegative constants $\mathbf
{c}=(\mathbf{c}_{\mathit{ii}'}, 1\leq i\neq i'\leq k)$ such that the
restriction of $\chi$ to $\{M\in[k]^{\mathbb{N}\otimes k}\dvtx
|M|_k=I_k\}$ is a
single-index flip measure,
\[
\mathbf{1}_{\{|M|_k=I_k\}}\chi(dM)=\sum_{1\leq i\neq i'\leq
k}
\mathbf {c}_{\mathit{ii}'}\rho_{\mathit{ii}'}.
\]
\end{longlist}

For (i), we let $\chi$ be the exchangeable characteristic measure of
$\mathbf{X}$ from Theorem~\ref{thmm:characteristic measure Feller}.
Then $\chi$ satisfies \eqref{eq:regularity chi} and we can write
$\chi
_n$ to denote the restriction of $\chi$ to the event $\{M\in
[k]^{\mathbb{N}\otimes k}
\dvtx M^{[n]}\neq\idnk\}$, for each $n\in\Nb$.
By \eqref{eq:regularity chi}, each $\chi_n$ is a finite measure on
$[k]^{\mathbb{N}\otimes k}$ and, by coset exchangeability, it is
invariant under action
by $k$-tuples of permutations $\sigma=(\sigma_1,\ldots,\sigma
_k)\dvtx\mathbb
{N}^k\rightarrow\mathbb{N}^k$ that fix $[n]^k$.
As a result, we define the $n$-shift $\overleftarrow{M}_{[n]}$ of
$M\in
[k]^{\mathbb{N}\otimes k}$ as follows: for $M:=(M_1,\ldots,M_k)$, we put
$\overleftarrow
{M}_{[n]}:=(\overleftarrow{M}_{1,[n]},\ldots,\overleftarrow
{M}_{k,[n]})$, where
\[
\overleftarrow{M}_{i,[n]}:=M_i^{n+1}M_i^{n+2}
\cdots,\qquad  i=1,\ldots,k.
\]
(The $n$-shift of $M$ is the coset decomposition of
$M'=M^{nk+1}M^{nk+2}\cdots,$ the $k$-coloring obtained by removing the
first $nk$ coordinates of $M$.)
The image $\overleftarrow{\chi}_n$ of $\chi_n$ by the $n$-shift is a
finite, coset exchangeable measure on $[k]^{\mathbb{N}\otimes k}$ that
satisfies~\eqref
{eq:regularity chi}.

By corollary to Theorem~\ref{thmm:cut-and-paste chain},
$\overleftarrow
{\chi}_n$-almost every $M\in[k]^{\mathbb{N}\otimes k}$ possesses
asymptotic frequency
$|M|_k\in\stochk$.
Since the asymptotic frequency of any $M\in[k]^{\mathbb{N}\otimes k}$
depends only on
its $n$-shift, for every $n\in\mathbb{N}$, $\chi_n$-almost every
$M\in
[k]^{\mathbb{N}\otimes k}$ possesses asymptotic frequency and, by
Theorem~\ref
{thmm:cut-and-paste chain}, we may write
%
%
\begin{equation}
\label{eq:n-finite measure} \chi_n(dM)=\int_{\stochk}
\mu_S(dM)\chi_n\bigl(|M|_k\in dS\bigr).
\end{equation}
Since $\chi_n\uparrow\chi$ as $n\uparrow\infty$, the monotone
convergence theorem implies that $\chi$-almost every $M\in
[k]^{\mathbb{N}\otimes k}$
possesses asymptotic frequencies.

To establish (ii), we consider the event that $\{M\in[k]^{\mathbb
{N}\otimes k}\dvtx
\overleftarrow{M}_{[n]}^{[2]}\neq\mbox{id}_{k,2}\}$ under~$\chi_n$.
(Here, $\overleftarrow{M}_{[n]}^{[m]}$ denotes the restriction to
$[k]^{[m]\otimes k}$ of the $n$-shift of $\overleftarrow{M}_{[n]}$.)
We define the $n$-shift measure by
%
%
\begin{equation}
\label{eq:n-shift measure} \overleftarrow{\chi}_n(dM)=\int_{\stochk}
\mu_{S}(dM)\overleftarrow {\chi }_n\bigl(\vert M
\vert_k\in dS\bigr),
\end{equation}
from which, for every $S\in\stochk$,
\begin{eqnarray*}
\chi_n\bigl(\bigl\{\overleftarrow{M}_{[n]}^{[2]}
\neq\operatorname{id}_{k,2}\bigr\} \mid \vert M\vert_k=S
\bigr)&=&\overleftarrow{\chi}_n\bigl({M}^{[2]}\neq
\operatorname {id}_{k,2}\mid\vert M\vert_k=S\bigr)
\\
&=&\mu_S\bigl(\bigl\{M^{[2]}\neq\operatorname{id}_{k,2}
\bigr\}\bigr)
\\
&\geq&1-S_*^2
\\
&\geq&1-S_*.
\end{eqnarray*}
Writing $\Sigma_n(dS):=\mathbf{1}_{\{\vert M\vert_k\neq I_k\}}\vert
{\chi
}_n\vert_k(dS)$, we obtain the inequality
%
%
\begin{equation}
\label{eq:chi ineq} \chi_n\bigl(\bigl\{\overleftarrow{M}_{[n]}^{[2]}
\neq\operatorname{id}_{k,2}\bigr\} \bigr)\geq \int_{\stochk}(1-S_*)
\Sigma_n(dS).
\end{equation}
By definition of ${\chi}_n$ and $\Sigma_n$, $\Sigma_n$ increases to
$\mathbf{1}_{\{\vert M\vert_k\neq I_k\}}\vert{\chi}\vert_k=:\Sigma
$ as
$n\rightarrow\infty$, the right-hand side above converges to
\[
\int_{\stochk}(1-S_*)\Sigma(dS),
\]
and $\Sigma(\{I_k\})=0$. On the other hand, the left-hand side in
\eqref
{eq:chi ineq} satisfies
\[
\chi_n\bigl(\bigl\{\overleftarrow{M}_{[n]}^{[2]}
\neq\operatorname{id}_{k,2}\bigr\} \bigr)\leq \chi\bigl(\bigl\{
\overleftarrow{M}_{[n]}^{[2]}\neq\operatorname{id}_{k,2}
\bigr\} \bigr)=\chi \bigl(\bigl\{M^{[2]}\neq\operatorname{id}_{k,2}
\bigr\}\bigr)<\infty,
\]
by coset exchangeability and \eqref{eq:regularity chi}. We conclude that
\[
\int_{\stochk}(1-S_*)\chi\bigl(\vert M\vert_k\in dS\bigr)=\int
_{\stochk} (1-S_*)\Sigma(dS)\leq\chi\bigl(\bigl\{
\overleftarrow{M}_{[n]}^{[2]}\neq \operatorname{id}_{k,2}
\bigr\}\bigr)<\infty;
\]
and $\Sigma$ satisfies \eqref{eq:regularity Sigma2}.

Finally, we must establish $\mathbf{1}_{\{\vert M\vert_k\neq I_k\}
}\chi
=\mu_\Sigma$. Indeed, for every $n\in\mathbb{N}$ and fixed $M^*\neq
\idnk
$, the monotone convergence theorem implies
\[
\chi\bigl(\bigl\{M^{[n]}=M^*, \vert M\vert_k\neq
I_k\bigr\}\bigr)=\lim_{m\uparrow\infty
}\chi\bigl(\bigl
\{M^{[n]}=M^*, \overleftarrow{M}_{[n]}^{[m]}\neq
\operatorname {id}_{k,m}, \vert M\vert_k\neq I_k
\bigr\}\bigr).
\]
By coset exchangeability, we can write
\[
\chi\bigl(\bigl\{M^{[n]}=M^*, \overleftarrow{M}_{[n]}^{[m]}
\neq\operatorname {id}_{k,m}, \vert M\vert_k\neq
I_k\bigr\}\bigr)=\overleftarrow{\chi}_m\bigl(\bigl\{
M^{[n]}=M^*, \vert M\vert_k\neq I_k\bigr\}
\bigr),
\]
and \eqref{eq:n-shift measure} implies
\[
\overleftarrow{\chi}_m\bigl(\bigl\{M^{[n]}=M^*, \vert M
\vert_k\neq I_k\bigr\}\bigr)=\int_{\stochk}
\mu_S\bigl(\bigl\{M^{[n]}=M^*\bigr\}\bigr)\overleftarrow{
\chi}_m\bigl(\vert M\vert _k\in dS\bigr),
\]
which converges to
\[
\int_{\stochk} \mu_S\bigl(\bigl
\{M^{[n]}=M^*\bigr\}\bigr)\Sigma(dS)=\mu_{\Sigma}\bigl(\bigl\{ M\in
[k]^{\mathbb{N}\otimes k}\dvtx M^{[n]}=M^*\bigr\}\bigr).
\]
As $n$ was chosen arbitrarily and the restriction $\vert M\vert_k\neq
I_k$ forbids $M=\idk$, we conclude (ii).

To establish (iii), let $\chi^*$ be the restriction of $\chi$ to the
event $\{M\in[k]^{\mathbb{N}\otimes k}\dvtx\break  M^{[2]}\neq
\operatorname{id}_{k,2}, \vert
M\vert_k=I_k\}$.
By \eqref{eq:regularity chi} and corollary to Theorem~\ref
{thmm:cut-and-paste chain}, $\chi^*$ is finite and its image
$\overleftarrow{\chi}^*_n$ by the $n$-shift is coset exchangeable;
thus, $\overleftarrow{\chi}^*_n$-almost every $M\in[k]^{\mathbb
{N}\otimes k}$ has
asymptotic frequency $\vert M\vert_k=I_k$ and $\overleftarrow{\chi
}^*_n$ is proportional to the unit mass at $\idk$.
So, we may restrict our attention to the event $E:=\{M^{[2]}\neq
\operatorname{id}_{k,2}, \overleftarrow{M}_{[3]}=\idk\}$ consisting
of maps $[k]^{\mathbb{N}}\rightarrow[k]^{\mathbb{N}}$ that fix coordinates
$n\geq3$.

Any $M=(M_1,\ldots,M_k)\in E$ is specified by a $k$-tuple
$((j_{11},j_{12}),(j_{21},j_{22}),\ldots,\break  (j_{k1},j_{k2}))$,
that is, the $i$th coset of $M$ [as in \eqref{eq:coset}] is
%
%
\begin{equation}
\label{eq:color pairs}M_i=j_{i1}j_{i2}\mathit{iii}\cdots,\qquad i=1,
\ldots,k.
\end{equation}
With $I=((j_{11},j_{12}),(j_{21},j_{22}),\ldots,(j_{k1},j_{k2}))$, we
write $M_I\in[k]^{\mathbb{N}}$ to denote the map in \eqref{eq:color pairs}.
Let $K:=\{((j_{11},j_{12}),\ldots,(j_{k1},j_{k2}))\}$ be the set of all
$k$-tuples and $K^*:=K\setminus\{I^*\}$, where $I^*\in K$ is defined as
\[
I^*:=\bigl((1,1),(2,2),\ldots,(k,k)\bigr).
\]
Then $E:=\bigcup_{I\in K^*} M_{I}$, which includes all single-index
flip maps $\kappa_{\mathit{ii}'}^{(n)}$ for $n=1,2$.

Now, since $\overleftarrow{\chi}^*_n$ is proportional to the point mass
at $\idk$, $\chi^*$ is the sum
\[
\chi^*(\cdot)=\sum_{I\in K^*}c_I
\delta_{M_I}(\cdot),
\]
where $\delta_{M_I}(\cdot)$ is the Dirac point mass at $M_I$.
By exchangeability, the requirement $\chi(\{M\dvtx M^{[2]}\neq
\operatorname
{id}_{k,2}\})<\infty$ forces $c_I=0$ unless $M_I$ is a single-index
flip map.
By extension of the above argument, any $M\in[k]^{\mathbb{N}\otimes
k}$ for which
$|M|_k=I_k$ and $\mathbf{c}_M>0$ must be a single-index flip map;
otherwise, by exchangeability, each index changes states at an infinite
rate and the finite restrictions cannot have c\`adl\`ag paths. This
establishes (iii) and completes the proof.
\end{pf*}

\subsection{Projection into the simplex}
By exchangeability of $\Xb$, the asymptotic frequency $|X_t|$ exists
almost surely for any fixed $t\geq0$.
In discrete-time, this and countable additivity of probability measures
imply the almost sure existence of $|\Xb|=(|X_m|, m\geq0)$.
In continuous-time, however, $\Xb=(X_t, t\geq0)$ is uncountable and
the corresponding conclusion does not follow immediately.
Nevertheless, Theorem~\ref{thmm:Levy-Ito} harnesses the behavior of
$\Xb
$ to a fruitful outcome: $|\Xb|=(|X_t|, t\geq0)$ exists and is a
Feller process.

To show this, we work on the compact metric space $(\simplexk,\tilde
{d})$, where
\[
\tilde{d}\bigl(s,s'\bigr):=\frac{1}{2}\sum
_{j=1}^k\bigl|s_j-s'_j\bigr|,\qquad
s,s'\in \simplexk.
\]
Under this metric, any $S\in\stochk$ determines a Lipschitz continuous
map $\simplexk\rightarrow\simplexk$, that is, for all $D,D'\in
\simplexk
$ and any $S\in\stochk$,
\[
\tilde{d}\bigl(DS,D'S\bigr)\leq\tilde{d}\bigl(D,D'
\bigr).
\]
We further exploit an alternative description of $\Xb^*_{\Sigma,\mathbf
{c}}$ by an associated Markov process on $[k]^{\mathbb{N}\otimes k}$.

Let $\Mb$ be the Poisson point process with intensity $dt\otimes\chi
_{\Sigma,\mathbf{c}}$, as above.
For each $n\in\Nb$, we define $\mathbf{F}^{[n]}:=(F_t^{[n]}, t\geq0)$
on $\mathop{[k]^{[n]\otimes k}} $ by $F_0^{[n]}=\idnk$ and:
\begin{itemize}
\item if $t>0$ is an atom time of $\Mb$ for which $M_t^{[n]}\neq\idnk$,
we put $F_t^{[n]}=M_t^{[n]}(F_{t-}^{[n]})$,
\item otherwise, we put $F_t^{[n]}=F_{t-}^{[n]}$.
\end{itemize}
We define $\mathbf{F}$ as the limit of $(\mathbf{F}^{[n]}, n\in\Nb)$,
which is a coset exchangeable, consistent Markov process on
$[k]^{\mathbb{N}\otimes k}$.
By its construction, $\mathbf{F}$ is closely tied
to $\Xb_{\Sigma,\mathbf{c}}^*=(X_t^*, t\geq0)$ by the relations:
\begin{itemize}
\item$|F_0|_k=I_k$ and
\item$X^*_t=F_t(X_0^*)$ for all $t\geq0$.
\end{itemize}

\begin{pf*}{Proof of Theorem~\ref{thmm:induced simplex process}}
Let $(\mathcal{F}_t, t\geq0)$ denote the natural filtration of
$\mathbf
{X}$ and, independently of $(\mathcal{F}_t, t\geq0)$, let $\mathbf
{F}:=(F_t, t\geq0)$ be the process on $[k]^{\mathbb{N}\otimes k}$
constructed above.
By Theorem~\ref{thmm:Levy-Ito}, the conditional law of $X_{t+s}$ given
$\mathcal{F}_t$ is that of $F_s(X_t)$.
By \eqref{eq:regularity chi} and exchangeability of $X_0$, $X_t$
possesses asymptotic frequencies almost surely for every $t\geq0$.
In fact, $|X_t|$ exists simultaneously for all $t\geq0$ with
probability one.

From Theorem~\ref{thmm:Levy-Ito}, a version of $\Xb$ can be constructed
as $\Xb^*_{\Sigma,\mathbf{c}}=(X_t^*, t\geq0)$, whose discontinuities
are of Types-(I) and (II) in Section~\ref{section:intro:Levy-Ito}.
In the projection $|\mathbf{X}^*_{\Sigma,\mathbf{c}}|$, discontinuities
only occur at the times of Type-(I) discontinuities, of which there are
at most countably many.
In between jumps, the trajectory of $|\mathbf{X}^*_{\Sigma,\mathbf
{c}}|$ is deterministic and continuous in $\simplexk$.
As a result, $|\mathbf{X}^*_{\Sigma,\mathbf{c}}|$ exists and is c\`
adl\`ag almost surely.
By corollary to Theorem~\ref{thmm:induced simplex chain}, $\vert
X_{t+s}^*\vert\equalinlaw|F_s(X_t^*)|=|F_s|_k|X_t^*|$, given
$\mathcal{F}_t$.
Since permutation does not affect the asymptotic frequency of either
$F_s$ or $X_t^*$, $|\mathbf{X}^*_{\Sigma,\mathbf{c}}|$ has the
Markov property.

Lipschitz continuity of every $S\dvtx\simplexk\rightarrow\simplexk$,
$S\in
\stochk$, implies the Feller property.
By compactness of $\simplexk$, any continuous $g\dvtx\simplexk
\rightarrow
\mathbb{R}$ is uniformly continuous and, therefore, bounded.
By the dominated convergence theorem, continuity of the map defined by
$S\in\stochk$, and Theorem~\ref{thmm:Levy-Ito}, the maps $D\mapsto
\mathbf{P}_t g(D)$ are continuous for all $t>0$.
By \eqref{eq:regularity chi}, $F_t\rightarrow\idk$ in probability as
$t\downarrow0$; whence, $\vert F_t\vert_k\rightarrow I_k$ and
$|F_t(X^*_0)|=|F_t|_k|X^*_0|\rightarrow|X^*_0|$, both in probability as
$t\downarrow0$.
We conclude that $\lim_{t\downarrow0}\mathbf{P}_tg(D)=g(D)$ for every
continuous function $g\dvtx\simplexk\rightarrow\mathbb{R}$, from which
follows the Feller property.
\end{pf*}

\section{Homogeneous cut-and-paste processes}\label{section:partition chains}

Theorems \ref{thmm:cut-and-paste chain}--\ref{thmm:induced simplex
process} extend to partition-valued processes with minor modifications.
Let $\Pib=(\Pi_t, t\geq0)$ be a continuous-time exchangeable,
consistent Markov process on $\partitionsNk$.
Specifically, $\Pib$ is a Markov process such that
\begin{longlist}[(A)]
\item[(A)] $\Pib^{\sigma}=(\Pi_t^{\sigma}, t\geq0)$ is a version of
$\Pib$ for all $\sigma\in\symmetricN$ and
\item[(B)] $\Pib^{[n]}=(\Pi^{[n]}_t, t\geq0)$ is a Markov chain on
$\mathop{\mathcal{P}_{[n]:k}} $, for every $n=1,2,\ldots.$
\end{longlist}
By Proposition~\ref{prop:Feller equiv}, $\Pib$ is a Feller process, and
thus, its evolution is determined by the finite jump rates
%
%
\begin{eqnarray}
\label{eq:fidi rates partition} Q_n\bigl(\pi,\pi'\bigr):=\lim
_{t\downarrow0}\frac{1}{t}\Pbb\bigl\{\Pi _t^{[n]}=
\pi'\mid \Pi_0^{[n]}=\pi\bigr\},
\nonumber
\\[-8pt]
\\[-8pt]
\eqntext{\pi\neq
\pi'\in\mathop{\mathcal{P}_{[n]:k}}, \mbox{ for each
}n\in\Nb,}
\end{eqnarray}
which satisfy \eqref{eq:finite markov rates}, \eqref{eq:exch markov
rates} and \eqref{eq:consistent markov rates}.

For any $\pi\in\partitionsNk$, we obtain its \emph{symmetric associate}
$\tilde{x}\in[k]^{\mathbb{N}}$ by labeling the blocks of $\pi$
uniformly and
without replacement in $[k]$.
In particular, for $\pi=(B_1,\ldots,B_r)\in\partitionsNk$ (listed in
order of least element), $\tilde{x}$ is a random $k$-coloring of $\Nb$
obtained by drawing labels $(l_1,\ldots,l_r)$ without replacement from
$[k]$ and putting $\tilde{x}=\tilde{x}^1\tilde{x}^2\cdots,$ where
\[
\tilde{x}^j=l_i \quad\Longleftrightarrow\quad j\in
B_i.
\]
Thus, $\Bb(\tilde{x})=\pi$ with probability one and each element in the
set $\mathcal{B}^{-1}(\pi)$ has equal probability.
For each $n\in\Nb$, we define the \emph{symmetric associate transition
rate} $\tilde{Q}_n$ on $\mathop{[k]^{[n]}} $ by
%
%
\begin{equation}
\label{eq:symmetric rates} \tilde{Q}_n\bigl(x,x'
\bigr):=Q_n\bigl(\Bb_n(x),\Bb_n
\bigl(x'\bigr)\bigr)/k^{\downarrow\#\Bb
_n(x')},\qquad x\neq x'\in
\mathop{[k]^{[n]}},
\end{equation}
where $\#\pi$ denotes the number of blocks of $\pi\in\partitionsN$ and
$k^{\downarrow j}:=k(k-1)\cdots(k-j+1)$.
Under $\tilde{Q}_n$, a transition from $x\in\mathop{[k]^{[n]}} $ is
obtained by
projecting $x\mapsto\Bb_n(x)=\pi$, generating a transition $\Pi
'\sim
Q_n(\pi,\cdot)$, and randomly coloring the blocks of $\Pi'$ to
obtain a
symmetric associate $\tilde{X}'\in\mathop{[k]^{[n]}} $.
The next proposition follows from definition \eqref{eq:symmetric rates}
and properties \eqref{eq:finite markov rates}, \eqref{eq:exch markov
rates} and \eqref{eq:consistent markov rates} of $(Q_n, n\in\Nb)$ in
\eqref{eq:fidi rates partition}.

%
\begin{prop}\label{prop:symmetric rates}
The collection $(\tilde{Q}_n, n\in\Nb)$ defined in \eqref{eq:symmetric
rates} determines a unique exchangeable transition rate measure $\tilde
{Q}$ on $[k]^{\mathbb{N}}$.
\end{prop}

From $\tilde{Q}$, we construct $\tilde{\Xb}=(\tilde{X}_t, t\geq0)$,
the \emph{symmetric associate} of $\Pib$, by first generating $\tilde
{X}_0$ as the symmetric associate of a partition from the initial
distribution of $\Pib$ and, given $\tilde{X}_0$, letting $\tilde{\Xb}$
evolve as a Markov process with initial state $\tilde{X}_0$ and
transition rate measure $\tilde{Q}$.

%
\begin{prop}\label{prop:symmetric associate}
The symmetric associate $\tilde{\Xb}$ of $\Pib$ is an exchangeable,
consistent Markov process on $[k]^{\mathbb{N}}$ and $\Bb(\tilde{\Xb
})=(\Bb
(\tilde{X}_t), t\geq0)$ is a version of~$\Pib$.
\end{prop}

\begin{pf}
We have constructed $\tilde{\Xb}$ so that it projects to and respects
the structure of $\Pib$.
To wit, $\Pib$ is exchangeable and consistent, and so is $\tilde{\Xb}$.
\end{pf}

For any permutation $\gamma\dvtx[k]\rightarrow[k]$, we define the
\emph{recoloring} of $x\in[k]^{\mathbb{N}}$ by
%
%
\begin{equation}
\label{eq:recoloring} \gamma x:=\gamma\bigl(x^1\bigr)\gamma
\bigl(x^2\bigr)\cdots.
\end{equation}
Since $\Bb(x)$ is the projection of $x$ into $\partitionsNk$ by
removing colors, recoloring does not affect $x\mapsto\Bb(x)$, that
is, $\Bb(x)=\Bb(\gamma x)$ for all $x\in[k]^{\mathbb{N}}$ and
$\gamma\in
\symmetrick$.
Thus, by definition \eqref{eq:symmetric rates}, $\tilde{Q}$ is
invariant under arbitrary recoloring of its arguments,
%
%
\begin{equation}
\label{eq:recoloring invariant} \tilde{Q}\bigl(\gamma x,\gamma'A\bigr)=
\tilde{Q}(x,A), \qquad x\in[k]^{\mathbb{N}}, A\subseteq[k]^{\mathbb{N}},
\end{equation}
for all $\gamma,\gamma'\in\symmetrick$, where $\gamma'A:=\{\gamma
'x'\dvtx x'\in A\}$ is the image of $A$ under recoloring by $\gamma'$.
By Theorem~\ref{thmm:characteristic measure Feller}, $\tilde{Q}$ is
characterized by a coset exchangeable measure $\tilde{\chi}$ which, by
condition \eqref{eq:recoloring invariant}, is invariant under the
action of \emph{left- and right-recoloring}, which we now define.

For $M\in[k]^{\mathbb{N}\otimes k}$ and $\gamma,\gamma'\in
\symmetrick$, we define the
\emph{left--right recoloring of $M$ by $(\gamma,\gamma')$} by
$M':=\gamma
M\gamma'$, where
%
%
\begin{equation}
\label{eq:row-column action} M'(x):=\gamma' M\bigl(
\gamma^{-1}x\bigr), \qquad x\in[k]^{\mathbb{N}},
\end{equation}
the $k$-coloring obtained by first recoloring $x$ by $\gamma^{-1}$,
then applying $M$, and finally recoloring by $\gamma'$.
We call a coset exchangeable measure \emph{row--column exchangeable} if
it is invariant under left--right recoloring by all pairs $(\gamma,\gamma
')\in\symmetrick\times\symmetrick$.

%
\begin{lemma}\label{lemma:RCE}
Let $\tilde{\chi}$ be the coset exchangeable measure that determines
$\tilde{Q}$.
Then $\tilde{\chi}$ is row--column exchangeable.
\end{lemma}

\begin{pf}
Fix $x\in[k]^{\mathbb{N}}$ and $A\subseteq[k]^{\mathbb{N}}$.
By \eqref{eq:recoloring invariant} and Theorem~\ref{thmm:Levy-Ito},
\begin{eqnarray*}
\tilde{\chi}\bigl(\bigl\{M\dvtx M(x)\in A\bigr\}\bigr)&=&\tilde{Q}(x,A)
\\
&=&\tilde{Q}\bigl(\gamma x,\gamma' A\bigr)
\\
&=&\tilde{\chi}\bigl(\bigl\{M\dvtx M(\gamma x)\in\gamma' A\bigr\}
\bigr)
\\
&=&\tilde{\chi}\bigl(\bigl\{M\dvtx\gamma^{-1}M\gamma'^{-1}(x)
\in A\bigr\}\bigr)
\\
&=&\tilde{\chi}\bigl(\bigl\{\gamma M\gamma'\dvtx M(x)\in A\bigr\}
\bigr),
\end{eqnarray*}
implying $\tilde{\chi}$ is row--column exchangeable.
\end{pf}

As a corollary to Theorem~\ref{thmm:Levy-Ito} and Proposition~\ref
{prop:symmetric associate}, $\tilde{\chi}$ is determined by a unique
pair $(\tilde{\Sigma},\tilde{\mathbf{c}})$, where $\tilde{\Sigma
}$ is a
measure satisfying \eqref{eq:regularity Sigma} and $\tilde{\mathbf
{c}}=(\tilde{\mathbf{c}}_{\mathit{ii}'}, 1\leq i\neq i'\leq k)$ is a collection
of nonnegative constants, that is,
%
%
\begin{equation}
\label{eq:tilde-chi} \tilde{\chi}=\mu_{\tilde{\Sigma}}+\sum
_{1\leq i\neq i'\leq
k}\tilde {\mathbf{c}}_{\mathit{ii}'}\rho_{\mathit{ii}'}.
\end{equation}
On $\stochk$, we call a measure $\Sigma$ \emph{row--column exchangeable}
if it is invariant under arbitrary permutation of rows and columns,
$S\mapsto\gamma S\gamma'^{-1}:=(S_{\gamma(i)\gamma'(i')}, 1 \leq
i,i'\leq k)$ for all $\gamma,\gamma'\in\symmetrick$.

%
\begin{prop}\label{prop:Sigma-RCE}
Let $\tilde{\chi}$ be as defined in \eqref{eq:tilde-chi}. Then
$\tilde
{\Sigma}$ is row--column exchangeable and there exists a unique $c\geq0$
such that $\tilde{\mathbf{c}}_{\mathit{ii}'}=c$ for all $1\leq i\neq i'\leq k$.
\end{prop}

\begin{pf}
In \eqref{eq:tilde-chi}, $\tilde{\chi}$ is expressed as the sum of
mutually singular measures, and we treat $\sum_{1\leq i\neq i'\leq
k}\tilde{\mathbf{c}}_{\mathit{ii}'}\rho_{\mathit{ii}'}$ first.

For $1\leq i\neq i'\leq k$ and $n\in\Nb$, we define
\[
A_{\mathit{ii}'}(n):=\bigl\{\kappa_{\mathit{ii}'}^{(1)},\ldots,
\kappa_{\mathit{ii}'}^{(n)}\bigr\},
\]
the subset of $[k]^{\mathbb{N}\otimes k}$ containing all single-index
flips from $i$ to
$i'$ for indices in $[n]$.
By Lemma~\ref{lemma:RCE}, $\tilde{\chi}$ is invariant under arbitrary
left- and right-recoloring as in \eqref{eq:row-column action}; whence,
\[
n\tilde{\mathbf{c}}_{\mathit{ii}'}=\tilde{\chi}\bigl(A_{\mathit{ii}'}(n)\bigr)=
\tilde{\chi }\bigl(A_{\gamma(i)\gamma(i')}(n)\bigr)=n\tilde{\mathbf{c}}_{\gamma(i)\gamma(i')}
\]
for all $n\in\Nb$ and $\gamma\in\symmetrick$, implying $\tilde
{\mathbf
{c}}_{\mathit{ii}'}=\tilde{\mathbf{c}}_{jj'}=c$ for all $i\neq i'$ and $j\neq j'$.

Restricted to the event $\{M\in[k]^{\mathbb{N}\otimes k}\dvtx
|M|_k\neq I_k\}$, $\tilde
{\chi}$ induces a measure $\tilde{\Sigma}$ satisfying \eqref
{eq:regularity Sigma} through the map $M\mapsto|M|_k$.
Row--column exchangeability follows by row--column exchangeability of
$\tilde{\chi}$ and definition of $M\mapsto|M|_k$ in \eqref
{eq:asymptotic frequency matrix}.
\end{pf}

\begin{pf*}{Proof of Theorem~\ref{thmm:Levy-Ito partition}}
For $\Pib$ in continuous-time, Theorem~\ref{thmm:Levy-Ito partition} is
a corollary of Theorem~\ref{thmm:Levy-Ito} and Propositions \ref
{prop:symmetric rates}, \ref{prop:symmetric associate} and \ref
{prop:Sigma-RCE}.
The discrete-time conclusion follows since single-index flips are not
permitted (forcing $c=0$) and Markov processes with finite jump rates
can be treated as discrete-time chains with exponentially distributed
hold times between jumps.
\end{pf*}

According to Theorem~\ref{thmm:induced simplex process}, the projection
into $\simplexk$ of an exchangeable $[k]^{\mathbb{N}}$-valued Feller process
exists and is also a Feller process.
The analogous projection of $\Pib$ into $\Delta_k^{\downarrow}$ by
$|\cdot
|^{\downarrow}$ also exists and is Feller.

\begin{pf*}{Proof of Theorem~\ref{thmm:induced simplex process-partition}}
Almost sure existence of $\vert\bolds{\Pi}\vert^{\downarrow}$ follows
from Theorem~\ref{thmm:Levy-Ito partition} and the existence of $\vert
\mathbf{X}\vert$ for any exchangeable Feller process on $[k]^{\mathbb{N}}$
(Theorem~\ref{thmm:induced simplex process}).
By Proposition~\ref{prop:Sigma-RCE}, the characteristic measure $\chi$
induces a row--column exchangeable measure $|\chi|_k$ on $\stochk$, and
so $\vert\bolds{\Pi}\vert^{\downarrow}$ is Markovian.
Theorem~\ref{thmm:induced simplex process} implies the Feller property
since $\vert\mathbf{X}\vert$ is Feller and any continuous
$g\dvtx\Delta_k^{\downarrow}\rightarrow\mathbb{R}$ induces a
continuous function
$g'\dvtx\simplexk\rightarrow\mathbb{R}$ which is symmetric in its arguments.
\end{pf*}

By the description in Theorem~\ref{thmm:Levy-Ito partition}, $\Pib$ is
characterized by its symmetric associate $\tilde{\Xb}$, whose
transition law treats colors homogeneously.
We commingle terms and call both $\tilde{\Xb}$ and $\Pib$ a \emph
{homogeneous cut-and-paste process} with parameter $(\tilde{\Sigma
},\tilde{c})$.

\subsection{Self-similar cut-and-paste processes}\label{section:Poisson
structure}

In \cite{Crane2011a}, we introduced a family of cut-and-paste chains,
which we now call self-similar homogeneous cut-and-paste chains.
We showed an instance of these chains in Example~\ref{ex:discrete}.

For a self-similar cut-and-paste process, the measure $\Sigma$ is the
$k$-fold product of some $\sigma$-finite measure on $\simplexk$, that
is, $\Sigma=\nu\otimes\cdots\otimes\nu$, for $\nu$ symmetric and
satisfying
%
%
\begin{equation}
\label{eq:regularity nu}\nu\bigl(\bigl\{(1,0,\ldots,0)\bigr\}\bigr)=0\quad
 \mbox{and}\quad \int
_{\Delta_k^{\downarrow}}(1-s_*)\nu(ds)<\infty,
\end{equation}
where $s_*:=\min\{s_1,\ldots,s_k\}$.
By symmetry of $\nu$, $\Sigma$ is row--column exchangeable.

The processes studied in \cite{Crane2011a} were \emph{pure-jump} in that
they did not admit single-index flips.
By letting single-index flips occur at rate $c\geq0$, we obtain the
class of \emph{self-similar homogeneous cut-and-paste processes} with
characteristic measure
\[
\chi=\mu_{\nu\otimes\cdots\otimes\nu}+c\rho,
\]
where $\rho:=\sum_{1\leq i\neq i'\leq k}\rho_{\mathit{ii}'}$.
The special case $c=0$ and $\nu=\PD(-\alpha/k,\alpha)$ plays a role in
clustering applications \cite{Crane2014a}.

\section{Concluding remarks}\label{section:conclusion}

\subsection{Relation to exchangeable coalescent and fragmentation
processes}\label{section:relation to cfp}
In spirit, our main theorems resemble previous results for exchangeable
coalescent and fragmentation processes.
In substance, our processes differ in fundamental ways.

\subsubsection{Bounded number of blocks}

All processes studied in this paper evolve on either $[k]^{\mathbb
{N}}$ or
$\partitionsNk$ for fixed $k\in\mathbb{N}$.
Bounding the number of blocks is necessary to characterize the jump
probabilities/rates by a measure on stochastic matrices.
Without an upper bound on the number of blocks, an exchangeable
partition need not admit proper asymptotic frequencies.
In general, for $\pi=\{B_1,B_2,\ldots\}\in\partitionsN$, the sum of its
asymptotic block frequencies may be strictly less than one, in which
case, it is common to write $s_0:=1-\sum_i |B_i|$ to denote the amount
of \emph{dust} in $|\pi|^{\downarrow}$.
For an exchangeable partition of $\mathbb{N}$, the dust is the totality
of its singleton blocks.
Furthermore, Theorem~\ref{thmm:Levy-Ito partition} requires the
cut-and-paste measure $\Sigma$ to treat all blocks symmetrically.
Without a uniform distribution on a countable set, we cannot specify
such a measure on $[k]^{\mathbb{N}\otimes k}$ with $k$ unbounded.

\subsubsection{Coalescent processes with finite initial state}

The representation in \eqref{eq:tilde-chi} covers a special subclass of
exchangeable coalescent processes whose initial state has a finite
number of blocks.
In this case, we let $k$ be the number of blocks of the initial state
$\Pi_0$, $c=0$, and $\Sigma$ a $\sigma$-finite row--column exchangeable
measure concentrated on $\{0,1\}$-valued stochastic matrices.
In this case, the homogeneous cut-and-paste process with initial state
$\Pi_0$ and characteristic measure $\chi=\mu_{\Sigma}$ is an
exchangeable coalescent.

On the other hand, no class of fragmentation processes corresponds to a
cut-and-paste process.
Fragmentation processes eventually fragment into the state of all
singletons, for which the number of blocks is infinite.

\subsubsection{Poissonian structure, coset mappings and $\Coag
$--$\Frag$
operators}

Exchangeable coalescent and fragmentation processes admit Poisson point
process constructions akin to our construction of $\Xb$ from the
Poisson point process $\Mb$ on $\mathbb{R}_+\times[k]^{\mathbb
{N}\otimes k}$.
For a coalescent process, $\mathbf{B}=\{(t,B_t)\}$ is a random subset
of $\mathbb{R}_+\times\partitionsN$ and $\Pib=(\Pi_t, t\geq0)$ is
constructed (informally) by putting $\Pi_t=\Coag(\Pi_{t-},B_t)$, for
each atom time $t$.
For $\pi,\pi'\in\partitionsN$,\break  $\Coag(\pi,\pi')$ is the \emph
{coagulation of $\pi$ by $\pi'$}, which determines a Lipschitz
continuous mapping $\partitionsN\rightarrow\partitionsN$.
Fragmentation processes have a similar construction in terms of the
$\Frag$-operator, which is also Lipschitz continuous.

The coset mappings, essential to our construction of cut-and-paste
processes, are also Lipschitz continuous.
To mimic the above constructions by the $\Coag$ and $\Frag$ operators,
we can define an operation $\cutpaste\dvtx\break [k]^{\mathbb{N}\otimes
k}\times\partitionsNk
\rightarrow\partitionsNk$ by
\[
\cutpaste(M,\pi):=\Bb\bigl(M(\tilde{x})\bigr),\qquad \tilde{x}
\mbox{ the symmetric associate of }\pi.
\]
From a Poisson point process $\Mb$ with intensity $dt\otimes\tilde
{\chi
}$, we generate $\Pib=(\Pi_t, t\geq0)$ (informally) by putting $\Pi
_t=\cutpaste(M_t,\Pi_{t-})$, for each atom time of $\Mb$.
The $\cutpaste$ operator differs from $\Coag$ and $\Frag$ because it
maps $[k]^{\mathbb{N}\otimes k}\times\partitionsNk\rightarrow
\partitionsNk$, rather
than $\partitionsN\times\partitionsN\rightarrow\partitionsN$.

We spare the details. See \cite{Crane2012c} for more on the interplay
between Poissonian structure, the Feller property and Lipschitz
continuous mappings.

\subsection{Equilibrium measures of cut-and-paste processes}

The process in Example~\ref{ex:discrete} is a self-similar homogeneous
cut-and-paste chain which is also reversible with respect to the
Poisson--Dirichlet distribution.
The process in Example~\ref{ex:continuous} evolves in continuous-time
and converges to a distribution whose projection to the simplex is
degenerate at $(1/2,1/2)$.
By Kingman's paintbox correspondence, these are the only possibilities.
In particular, the unique equilibrium measure of an exchangeable
cut-and-paste process, if it exists, is one of Kingman's paintbox measures.
The cut-and-paste representation is a powerful tool for studying
equilibrium measures of these chains, evinced by Crane and Lalley \cite
{CraneLalley2012a}.

\section*{Acknowledgements}

I thank Steve Lalley for his helpful comments on an early version and
Chris Burdzy for his kind encouragement.

%



\printaddresses


\begin{thebibliography}{16}

%
%
\bibitem{AldousExchangeability}
\begin{bincollection}[mr]
\bauthor{\bsnm{Aldous},~\bfnm{David~J.}\binits{D.~J.}}
(\byear{1985}).
\btitle{Exchangeability and related topics}.
In \bbooktitle{\'{E}cole D'\'et\'e de Probabilit\'es de
{S}aint-{F}lour, {XIII}---1983}.
\bseries{Lecture Notes in Math.}
\bvolume{1117}
\bpages{1--198}.
\bpublisher{Springer},
\blocation{Berlin}.
\bid{mr={0883646}}
\end{bincollection}
\bptok{imsref}
\endbibitem

%
%
\bibitem{Bertoin2001a}
\begin{barticle}[mr]
\bauthor{\bsnm{Bertoin},~\bfnm{Jean}\binits{J.}}
(\byear{2001}).
\btitle{Homogeneous fragmentation processes}.
\bjournal{Probab. Theory Related Fields}
\bvolume{121}
\bpages{301--318}.
\bid{issn={0178-8051}, mr={1867425}}
\end{barticle}
\bptok{imsref}
\endbibitem

%
%
\bibitem{Bertoin2006}
\begin{bbook}[mr]
\bauthor{\bsnm{Bertoin},~\bfnm{Jean}\binits{J.}}
(\byear{2006}).
\btitle{Random Fragmentation and Coagulation Processes}.
\bseries{Cambridge Studies in Advanced Mathematics}
\bvolume{102}.
\bpublisher{Cambridge Univ. Press},
\blocation{Cambridge}.
\bid{mr={2253162}}
\end{bbook}
\bptok{imsref}
\endbibitem

%
%
\bibitem{Crane2014a}
\begin{bmisc}[auto:STB|2014/02/12|14:17:21]
\bauthor{\bsnm{Crane},~\bfnm{H.}\binits{H.}}
(\byear{2014}).
\bhowpublished{Clustering from partition data. Unpublished manuscript.}
\end{bmisc}
\bptok{imsref}
\endbibitem

%
%
\bibitem{Crane2012c}
\begin{bmisc}[auto:STB|2014/02/12|14:17:21]
\bauthor{\bsnm{Crane},~\bfnm{H.}\binits{H.}}
(\byear{2014}).
\bhowpublished{Lipschitz partition processes. \emph{Bernoulli}. To appear.}
\end{bmisc}
\bptok{imsref}
\endbibitem

%
%
\bibitem{Crane2011a}
\begin{barticle}[mr]
\bauthor{\bsnm{Crane},~\bfnm{Harry}\binits{H.}}
(\byear{2011}).
\btitle{A consistent {M}arkov partition process generated from the
paintbox process}.
\bjournal{J. Appl. Probab.}
\bvolume{48}
\bpages{778--791}.
\bid{issn={0021-9002}, mr={2884815}}
\end{barticle}
\bptok{imsref}
\endbibitem

%
%
\bibitem{CraneLalley2012a}
\begin{barticle}[mr]
\bauthor{\bsnm{Crane},~\bfnm{Harry}\binits{H.}} \AND
\bauthor{\bsnm{Lalley},~\bfnm{Steven~P.}\binits{S.~P.}}
(\byear{2013}).
\btitle{Convergence rates of {M}arkov chains on spaces of partitions}.
\bjournal{Electron. J. Probab.}
\bvolume{18}
\bpages{23}.
\bid{issn={1083-6489},mr={3078020} }
\end{barticle}
\bptok{imsref}
\endbibitem

%
%
\bibitem{DurrettGranovskyGueron1999}
\begin{barticle}[mr]
\bauthor{\bsnm{Durrett},~\bfnm{Richard}\binits{R.}},
\bauthor{\bsnm{Granovsky},~\bfnm{Boris~L.}\binits{B.~L.}} \AND
\bauthor{\bsnm{Gueron},~\bfnm{Shay}\binits{S.}}
(\byear{1999}).
\btitle{The equilibrium behavior of reversible
coagulation-fragmentation processes}.
\bjournal{J. Theoret. Probab.}
\bvolume{12}
\bpages{447--474}.
\bid{issn={0894-9840}, mr={1684753}}
\end{barticle}
\bptok{imsref}
\endbibitem

%
%
\bibitem{Ewens1972}
\begin{barticle}[mr]
\bauthor{\bsnm{Ewens},~\bfnm{W.~J.}\binits{W.~J.}}
(\byear{1972}).
\btitle{The sampling theory of selectively neutral alleles}.
\bjournal{Theoret. Population Biology}
\bvolume{3}
\bpages{87--112};
\bnote{erratum, ibid. \textbf{3} (1972), 240; erratum, ibid.
\textbf{3} (1972), 376}.
\bid{issn={0040-5809}, mr={0325177}}
\end{barticle}
\bptok{imsref}%
\endbibitem

%
%
\bibitem{FurstenbergKesten1960}
\begin{barticle}[mr]
\bauthor{\bsnm{Furstenberg},~\bfnm{H.}\binits{H.}} \AND
\bauthor{\bsnm{Kesten},~\bfnm{H.}\binits{H.}}
(\byear{1960}).
\btitle{Products of random matrices}.
\bjournal{Ann. Math. Statist.}
\bvolume{31}
\bpages{457--469}.
\bid{issn={0003-4851}, mr={0121828}}
\end{barticle}
\bptok{imsref}
\endbibitem

%
%
\bibitem{Kingman1978a}
\begin{barticle}[mr]
\bauthor{\bsnm{Kingman},~\bfnm{J.~F.~C.}\binits{J.~F.~C.}}
(\byear{1978}).
\btitle{Random partitions in population genetics}.
\bjournal{Proc. R. Soc. Lond. Ser. A Math. Phys. Eng. Sci.}
\bvolume{361}
\bpages{1--20}.
\bid{issn={0962-8444}, mr={0526801}}
\end{barticle}
\bptok{imsref}
\endbibitem

%
%
\bibitem{Kingman1978b}
\begin{barticle}[mr]
\bauthor{\bsnm{Kingman},~\bfnm{J.~F.~C.}\binits{J.~F.~C.}}
(\byear{1978}).
\btitle{The representation of partition structures}.
\bjournal{J. Lond. Math. Soc. (2)}
\bvolume{18}
\bpages{374--380}.
\bid{issn={0024-6107}, mr={0509954}}
\end{barticle}
\bptok{imsref}
\endbibitem

%
%
\bibitem{Kingman1982}
\begin{barticle}[mr]
\bauthor{\bsnm{Kingman},~\bfnm{J.~F.~C.}\binits{J.~F.~C.}}
(\byear{1982}).
\btitle{The coalescent}.
\bjournal{Stochastic Process. Appl.}
\bvolume{13}
\bpages{235--248}.
\bid{issn={0304-4149},mr={0671034} }
\end{barticle}
\bptok{imsref}
\endbibitem

%
%
\bibitem{Pitman1999a}
\begin{barticle}[mr]
\bauthor{\bsnm{Pitman},~\bfnm{Jim}\binits{J.}}
(\byear{1999}).
\btitle{Coalescents with multiple collisions}.
\bjournal{Ann. Probab.}
\bvolume{27}
\bpages{1870--1902}.
\bid{issn={0091-1798}, mr={1742892}}
\end{barticle}
\bptok{imsref}
\endbibitem

%
%
\bibitem{Pitman2002a}
\begin{barticle}[mr]
\bauthor{\bsnm{Pitman},~\bfnm{Jim}\binits{J.}}
(\byear{2002}).
\btitle{Poisson--{D}irichlet and {GEM} invariant distributions for
split-and-merge transformation of an interval partition}.
\bjournal{Combin. Probab. Comput.}
\bvolume{11}
\bpages{501--514}.
\bid{issn={0963-5483}, mr={1930355}}
\end{barticle}
\bptok{imsref}
\endbibitem

%
%
\bibitem{Pitman2005}
\begin{bbook}[mr]
\bauthor{\bsnm{Pitman},~\bfnm{J.}\binits{J.}}
(\byear{2006}).
\btitle{Combinatorial Stochastic Processes}.
\bseries{Lecture Notes in Math.}
\bvolume{1875}.
\bpublisher{Springer},
\blocation{Berlin}.
\bid{mr={2245368}}
\end{bbook}
\bptok{imsref}
\endbibitem

\end{thebibliography}
\end{document}